\documentclass[twoside,12pt]{article}
\usepackage[usenames,dvipsnames]{color}
\usepackage{amsmath}
\usepackage{amsthm}
\usepackage{amssymb}
\usepackage[mathcal]{euscript}
\usepackage{mathrsfs}
\usepackage{verbatim}
\usepackage{cite}

\definecolor{viola}{rgb}{0.3,0,0.7}
\definecolor{ciclamino}{rgb}{0.5,0,0.5}

\def\pier #1{{\color{red}#1}}
\def\pier #1{#1}

\def\colli #1{{\color{blue}#1}}
\def\colli #1{#1}

\def\revis #1{{\color{red}#1}}
\def\revis #1{#1}

\def\sign{\mathop{\rm sign}}

\def\Sign{\mathop{\rm Sign}}

\newtheorem{Lemma1}{Lemma}[section]

\newtheorem{Teo-esistenza}[Lemma1]{Theorem - (Existence)}

\newtheorem{theorem}{Theorem}[section]

\title{Global existence for a singular phase field system \\
related to a sliding mode control problem\/\footnote{\pier{{\bf 
Acknowledgments.}\quad\rm The first author gratefully acknowledges some 
financial support from from the MIUR-PRIN Grant 2015PA5MP7 
``Calculus of Variations'', the GNAMPA (Gruppo Nazionale per l'Analisi 
Matematica, la Probabilit\`a e le loro Applicazioni) of INdAM (Istituto 
Nazionale di Alta Matematica) and the IMATI -- C.N.R. Pavia.}}}

\author{Pierluigi Colli\\
Dipartimento di Matematica, Universit\`a degli Studi di Pavia\\
Via Ferrata~5, 27100 Pavia, Italy\\
E-mail: \texttt{pierluigi.colli@unipv.it}\\
\and \\ Michele Colturato\\
Dipartimento di Matematica, Universit\`a degli Studi di Pavia\\
Via Ferrata~5, 27100 Pavia, Italy\\
E-mail: \texttt{michele.colturato01@universitadipavia.it}}
 
\date{}
\newcommand\testopari{\sc P. Colli and M. Colturato}
\newcommand\testodispari{\sc Singular system
related to a sliding mode control problem}

\begin{document}
\maketitle

\begin{abstract}
In the present contribution we consider a singular phase field system 
located in a \revis{smooth and bounded three-dimensional domain}.
The entropy balance equation is perturbed by 
a logarithmic nonlinearity \pier{and by the presence of an additional term involving a 
possibly nonlocal maximal monotone operator and arising from a class of sliding mode control problems.}
The second equation of the system accounts for the phase dynamics\revis{, and it} is deduced from
a balance law for the microscopic forces that are responsible for the phase transition process.
The resulting system is highly \revis{nonlinear; the} main difficulties lie in the \pier{contemporary 
presence of two nonlinearities, one \colli{of which} under time derivative, in the entropy balance equation}. Consequently, we are able to prove only the existence of \colli{solutions}.
To this aim, we will introduce a backward finite differences scheme and 
\pier{argue on this by proving uniform estimates and passing to the limit on the time step.}

\vspace{2mm}
\noindent \textbf{Key words:}~~ Phase field system; \pier{maximal monotone nonlinearities; nonlocal terms; initial and boundary value problem; existence of solutions}. 

\vspace{2mm}
\noindent \textbf{AMS (MOS) subject clas\-si\-fi\-ca\-tion:} 35K61, \pier{35K20,}
35D30, 80A22.

\end{abstract}

\section{Introduction}
This paper is devoted to the mathematical analysis of a system of partial differential equations (PDE) arising from a thermodynamic model describing phase
transitions. The system is written in terms of a rescaled balance of
energy and of a balance law for the microforces that govern the phase transition.
Moreover, the first equation of the system is perturbed by the presence of an additional maximal monotone nonlinearity. This paper will focus only on analytical aspects and,
in particular, will investigate the existence of solutions. 
In order to make the presentation clear from the beginning, we briefly introduce the main ingredients of the PDE system and give some comments on the physical meaning.

We \revis{deal with} a two-phase system located in a smooth bounded domain $\Omega \subseteq \mathbb{R}^3$
and let $T > 0$ denote \revis{some} final time. The unknowns of the problem are the 
\revis{\emph{absolute temperature} $\vartheta$ and an \emph{order parameter}} 
$\chi$ which \revis{can} represent the local proportion of one of the two phases. To
ensure thermomechanical consistency, suitable physical constraints on $\chi$ are \revis{considered}: if it is
assumed, e.g., that the two phases may coexist at each point with different proportions, it turns
out to be reasonable to require that $\chi$ lies between $0$ and $1$, with
$1 - \chi$ representing the proportion of the second phase. In particular, 
the values $\chi= 0$  and $\chi= 1$ may correspond to the pure phases, 
while $\chi$  is between $0$ and $1$ in the regions when both phases are
present. Clearly, the \revis{the system provides} an evolution for
$\chi$ that \revis{has to comply} with the previous physical constraint. 

Now, let us state precisely the equations as well as the initial and boundary conditions.
\revis{The equations} governing the evolution of $\vartheta$ and $\chi$ are recovered as balance laws. 
The first equation \revis{comes from} a reduction of the energy balance
equation divided by the absolute temperature $\vartheta$ (see \cite[formulas~(2.33)--(2.35)]{COLLIlog2}). 
\revis{Therefore}, the so-called entropy balance can be written in $\Omega \times (0, T )$ as follows:
\begin{equation} \label{II144}
\partial_t( \ln{\vartheta} + \ell \chi) - k_0 \Delta \vartheta = F,
\end{equation}
where \revis{$\ell$ is a positive parameter, $k_0 > 0$ is a thermal coefficient for the entropy flux $\textbf{Q}$, 
which is related to the heat flux vector $\textbf{q}$ by $\textbf{Q} = \textbf{q}/ \vartheta$,
 and $F$ stands for an external entropy source.}

In the present contribution, we assume that the entropy balance equation \eqref{II144} is perturbed by the presence of an additional maximal monotone nonlinearity, i.e., 
\begin{equation} \label{II1}
\partial_t( \ln{\vartheta} + \ell \chi) - k_0 \Delta \vartheta + \zeta= F,
\end{equation}
where
\begin{equation} \label{II3}
\zeta(t) \in A(\vartheta(t) - \vartheta^*) \quad \quad \textrm{for a.e. $t \in (0,T)$}. 
\end{equation}
Here, $\vartheta^*$ is a positive and smooth function ($\vartheta^*\in H^2(\Omega)$ with null outward normal derivative on the boundary) and  $A: L^2(\Omega) \to   L^2(\Omega)$ is a maximal monotone operator satisfying some conditions, namely:
$A$ is the subdifferential of a proper, convex and lower 
semicontinuous (l.s.c.) function $\Phi:  L^2(\Omega) \to  \mathbb{R}$ which
takes \colli{its} minimum in $0$, and $A$ is linearly bounded in $L^2(\Omega)$. 
In order to explain the role of this further nonlinearity,
we refer to \cite{BaCoGiMaRo}, where a class of sliding mode control problems is considered:
a state-feedback control $(\vartheta,\chi) \mapsto u(\vartheta,\chi)$ 
is added in the balance equations with the purpose of forcing the trajectories of the system 
to reach the sliding surface (i.e., a manifold of lower dimension where the control
goal is fulfilled and such that the original system restricted to this manifold has a
desired behavior) in finite time and maintains them on it.
As widely described in \cite{BaCoGiMaRo},
this study is physically meaningful in the framework of phase transition processes.

Let us mention the contributions \cite{Michele, Michele2}, where standard phase field 
systems of Caginalp type, 
perturbed by the presence of nonlinearities similar to \eqref{II3}, are considered. 
In \cite{Michele, Michele2} the existence of strong solutions, the global well-posedness 
of the system and the sliding mode property can be proved; unfortunately, here the problem we consider is rather more delicate due to the doubly nonlinear character of equation \eqref{II1} and it turns out that we cannot perform a so complete analysis. On the other hand,
we observe that, due to the presence of the logarithm of the temperature 
in the entropy equation \eqref{II1}, in the system \colli{we investigate} here 
the positivity of the variable representing the absolute temperature follows directly
from solving the problem, i.e., from finding a solution component $\vartheta$ to 
which the logarithm applies. This is \revis{an important feature and avoids} the use of other 
methods or the setting of special assumptions, in order to guarantee the positivity 
of  $\vartheta$ in the space-time domain. 

The second equation of the system under study \revis{describes} the phase dynamics 
and is deduced from a balance law for the microscopic forces that are 
responsible for the phase transition process.
According to \cite{fremond, Gurtin}, this balance reads
\begin{equation} \label{II2}
\partial_t \chi - \Delta \chi + \beta(\chi) + \pi(\chi) \ni \ell \vartheta,
\end{equation}
where $\beta + \pi$ represents the derivative, or the subdifferential,  of a double-well potential $\mathcal{W}$
defined as 
\begin{equation} \nonumber
\mathcal{W}= \tilde{\beta} + \tilde{\pi},
\end{equation}
where
	\begin{equation} \label{Intro-beta}
\tilde{\beta}: \mathbb{R} \longrightarrow [0, + \infty] \textrm{ is proper, l.s.c. and convex with $\tilde{\beta}(0)=0$,} 
		\end{equation}
	\begin{equation} \label{Intro-lambda}
	\textrm{$ \tilde{\pi} \in C^1(\mathbb{R})$ and $\pi= \tilde{\pi}' $ is Lipschitz continuous in $\mathbb{R}$.}  
	\end{equation}
{}\colli{Due to \eqref{Intro-beta},} the subdifferential $ \beta:= \partial \tilde{\beta}$ is well defined 
and \colli{turns out to be} a maximal 
monotone graph.	Moreover, as $\tilde{\beta}$ takes on its minimum in $0$, we have that 
$0 \in \beta(0)$. Note that in \eqref{II2} the inclusion is used in place of the equality in order to allow for the presence of a multivalued $\beta$. 

We recall that many different choices of $\tilde{\beta}$ and $\tilde{\pi}$ have been 
introduced in the literature (see, e.g., \cite{Blowey, COLLIlog3, Fabrizio, Penrose}). 
In \revis{case of} a solid-liquid phase transition, $\mathcal{W}$ may be \revis{taken in a} way 
that the full potential (cf.~\eqref{II2})
$$
\chi \mapsto \tilde{\beta} (\chi)  + \tilde{\pi} (\chi) - \ell \vartheta \chi 
$$
exhibits one of the two minima $ \chi= 0$ and $\chi = 1$ as global minimum for equilibrium, 
depending on whether $\vartheta $ is below or above a critical value $\vartheta_c$, which
may represent a phase change temperature. A sample case is given by
$\tilde{\pi} (\chi) =\ell \vartheta_c \chi$ and by the $\tilde{\beta}$ that coincides with the  indicator function $ I_{[0, 1]}$ of the interval $[0,1]$, that is, 
\begin{equation*}
\tilde{\beta} (\rho) =  I_{[0,1]} (\rho) =  \ \left\{
\begin{array}{ll}
0 & \text{if \ $0\leq \rho \leq 1$}
\\[0.1cm]
+\infty \ &\text{elsewhere}
\end{array}
\right. 
\end{equation*}
so that $\beta = \partial I_{[0, 1]}$ is specified by
\begin{equation*}
r \in \beta (\rho) 
\quad \hbox{ if and only if } \quad 
r \ \left\{
\begin{array}{ll}
\displaystyle
\leq \, 0 \   &\hbox{if } \ \rho=0   
\\[0.1cm]
= \, 0 \   &\hbox{if } \ 0< \rho < 1  
\\[0.1cm]
\geq \, 0 \  &\hbox{if } \  \rho  = 1  
\\[0.1cm]
\end{array}
\right. .
\end{equation*}
Of course, this yields a singular case for the potential $\colli{\cal W}$, in which $\tilde{\beta}$
is not differentiable, and it is known in the literature as the double obstacle case  (cf.~\cite{Blowey, COLLIlog3, fremond})

In the last decades phase field models have attracted a number of mathematicians and 
applied scientists to describe many different physical phenomena. Let us \colli{just} recall some results in the literature that are related to our system. Some key references are the papers \cite{COLLIlog, COLLIlog2, COLLIlog3}. Besides, we quote \cite{COLLIlog4}, 
where a first simplified version of the entropy system is considered,
and \cite{BCG,BFR} for related analyses and results. About special choices of the heat flux 
and phase field models ensuring positivity of the absolute temperature, we aim to quote the papers \cite{cl, cls, COLL-Sprek}, where some \revis{Penrose--Fife} models have been addressed.

The full problem investigated in this paper consists of equations \eqref{II1}--\eqref{II2} 
coupled with suitable boundary and initial conditions. In particular, we prescribe a no-flux condition on the boundary for both variables:
\begin{equation} \label{II5}
\partial_{\nu} \vartheta =0, \quad  \partial_{\nu}\chi = 0 \quad \textrm{ on $\Gamma \times (0, T )$},
\end{equation}
where $\partial_{\nu}$ denotes the outward normal derivative on the boundary $\Gamma$ of $\Omega$. Besides, in the light of \eqref{II3}, initial conditions are stated for $\ln \vartheta$ and $\chi\/$:
\begin{equation} \label{II6}
\ln{\vartheta}(0) = \ln{\vartheta_0}, \quad \chi(0) = \chi_0 \quad \textrm{in $\Omega$}.
\end{equation}

The resulting system is highly nonlinear\revis{. The} main difficulties lie in the treatment 
of the doubly nonlinear equation \eqref{II1}. The expert reader can realise that it is 
not trivial to recover some coerciveness and regularity for $\vartheta$ from \eqref{II1}, \eqref{II3} and 
\eqref{II5}; morever, the presence of both $\ln \vartheta$ under time derivative and the selection $\zeta $ from $A(\vartheta - \vartheta^*)$ complicates possible uniqueness arguments. 
For the moment, we \colli{are} just able to prove
the existence of solutions for the described problem.
To this aim, \colli{we introduce} a backward finite differences scheme and first examine the solvability of it, for which we have to introduce another approximating problem based on the use of Yosida regularizations for the maximal monotone operators.

As far as the outline of the paper is concerned, we state precisely assumptions and main results in Section~2, then introduce the time-discrete problem $(P_{\tau})$ in Section~3 and completely prove existence and uniqueness of the solution.
Section 4 is devoted to the proof of several uniform estimates, independent of $\tau$, involving the solution of $(P_{\tau})$. Finally, in Section 5 we pass to the limit as
as $\tau \searrow 0$ by means of compactness and monotonicity arguments in order to find a solution to the problem  \eqref{II1}--\eqref{II2}, \eqref{II5}--\eqref{II6}.

\section{Main results}
\setcounter{equation}{0}

\subsection{Preliminary assumptions}
We assume $\Omega \subseteq \mathbb{R}^3$ to be open, bounded, connected, of class $C^1$ and we write $| \Omega |$ for its Lebesgue measure.
Moreover, $\Gamma$ and $\partial_{\nu}$ stand for the boundary of $\Omega$ and the outward normal derivative, respectively.
Given a finite final time $T > 0$, for every $t \in (0,T]$ we set
\begin{equation*}
Q_t =  \Omega \times (0,t) , \ \ Q = Q_T,  \quad
\Sigma_t = \Gamma \times (0,t), \ \ \Sigma = \Sigma_T. 
\end{equation*}
We also \pier{introduce the spaces}
\begin{align}
H = L^2(\Omega), \quad V = H^1(\Omega), 
\label{W}
\quad W = \{ u \in H^2(\Omega): \ \partial_{\nu} u = 0 \ \textrm{on} \ \Gamma \},
\end{align}
with usual norms $\| \cdot \|_{H}$, $\|\cdot \|_{V}$, $\|\cdot \|_{W}$
and related inner products $(\cdot,\cdot )_{H}$, $(\cdot ,\cdot )_{V}$, $(\cdot ,\cdot )_{W}$, respectively. \pier{We identify $H$ with its dual space $H'$, so that 
$W\subset V \subset H \subset V' \subset W'$ with dense and compact embeddings. Let 
$\langle \cdot, \cdot \rangle$ \colli{denote} the duality pairing between $V'$ and $V$.} 
The notation $\| \cdot \|_p$ $(1 \leq p \leq \infty)$ stands for the standard norm
in $L^p(\Omega)$. For short, \pier{in the notation of norms} we do not distinguish between a space and a power thereof.

\pier{From now on, we interpret the operator $- \Delta$ as the Laplacian operator from the space $W$ to $H$, then including the Neumann homogeneous boundary condition. Moreover,
we extend $- \Delta$ to an operator from $V$ to $V'$ by setting 
\begin{equation} \label{pier4}
\langle -\Delta u, v \rangle := \int_\Omega \nabla u \cdot \nabla v, \quad u,v\in V.
\end{equation}%
}%
Throughout the paper, we account for the well-known \pier{continuous} embeddings $V \subset L^q(\Omega)$, with
$1 \leq q \leq 6$, $W \subset C^0(\overline{\Omega})$ and for the related Sobolev inequalities:
\begin{equation*}
\| v \|_q \leq \pier{C_s}\, \| v \|_V \quad \quad \textrm{and} \quad \quad \| v \|_{\infty} \leq \pier{C_s} \, \| v \|_W
\end{equation*}
for $v \in V$ and $v \in W$, respectively, where $\pier{C_s}$ depends on $\Omega$ only, since sharpness is not needed.
We will also use a variant of the Poincar\'e inequality, i.e., there exists a positive constant $\pier{C_p}$ such that
\begin{equation} \label{Poincare}
\| v \|_V \leq \pier{C_p} \Big( \| v \|_{L^1(\Omega)} + \| \nabla v \|_H \Big), \quad \quad v \in V.
\end{equation}
Furthermore, we make repeated use of the H\"older inequality, and of Young's inequalities, i.e.,
for every $a,b >0$, $\alpha \in (0,1)$ and $\delta > 0$ we have that 
\begin{gather} \label{dis1}
ab \leq \alpha a^{\frac{1}{\alpha}} + (1 - \alpha)b^{\frac{1}{1- \alpha}},
\\\label{dis2}
ab \leq \delta a^2 + \frac{1}{4\delta}b^2.
\end{gather}
\pier{Besides}, for every $a$, $b \in \mathbb{R}$ we have that
\begin{equation} \label{dis777} 
(a- b)a = \frac{1}{2}a^2 - \frac{1}{2}b^2 + \frac{1}{2}(a-b)^2  .
\end{equation}
We also recall the discrete version of the Gronwall lemma (see, e.g., \cite[Prop. 2.2.1]{Jerome}).
\begin{Lemma1}\label{grondiscr} 
If $(a_0, \colli{\ldots\,}, a_N) \in [0, + \infty)^{N+1}$ and $(b_1, \colli{\ldots\,}, b_N) \in [0, + \infty)^{N}$ satisfy 
\begin{equation*} 
a_m \leq a_0 + \sum_{n=1}^{m-1}  a_n b_n  \quad \quad \hbox{for \, $m = 1, \colli{\ldots\,}, N$,}
\end{equation*}
then
\begin{equation} \label{pier1} 
a_m \leq a_0 \ \textrm{exp} \Bigg( \sum_{n=1}^{m-1}   b_n \Bigg) \quad \quad \textrm{for \,$m = 1, \colli{\ldots\,}, N$.}
\end{equation}
\end{Lemma1}

\pier{Finally, we state another useful result for the sequel.}

\begin{Lemma1}\label{superlemma1}
Assume that $a$, $b \in \mathbb{R}$ are strictly positive. Then
\begin{equation} \label{pier2}
(a -b) \leq \big(\ln a^2 - \ln b^2\big) (a + b).
\end{equation}
\end{Lemma1}

\noindent
\textit{Proof}. We consider $a >b$ (if $b > a$ the technique of the proof is analogous) and obtain
\begin{equation*} \label{lemmm11}
(a -b) \leq \colli{{}(\ln a^2 - \ln b^2) (a + b) = 2( \ln a -  \ln b)} (a + b) = 2  \ln\bigg(\frac{a}{b} \bigg) (a + b).
\end{equation*}
\pier{Then, dividing} by $b$, we have that
\begin{equation} \label{lemmm12}
\bigg(  \frac{a}{b} - 1 \bigg) \leq 2  \ln\bigg( \frac{a}{b} \bigg) \bigg(  \frac{a}{b} + 1 \bigg).
\end{equation}
\pier{Letting} $x= a/b$, we can rewrite \eqref{lemmm12} as
\begin{equation*} \label{lemmm13}
( x - 1 ) \leq 2  \colli{ (  x + 1 ) \ln x }\quad \  \textrm{for $x \geq 1$.}
\end{equation*}
Now, we observe that \eqref{pier2} is verified if and only if the function
\begin{equation} \label{lemmm14}
f(x) := 2  (  x + 1 ) \ln x   - x + 1   \quad \textrm{\pier{is nonnegative} 
for every $x \geq 1$.}
\end{equation}
Since \pier{$f(1)=0$ and $f'(x) > 0$ for every $x \geq 1$},
we conclude that \eqref{lemmm14} holds. Then, the proof of the lemma is complete. $\Box$

In the following, the small-case symbol $c$ stands for different constants which depend only on $\Omega$, on the final time $T$,
on the shape of the nonlinearities and on the constants and the norms of the functions involved in the assumptions of our statements.
On the contrary, we use different symbols to denote precise constants to which we could refer.
\revis{It is important to point out} that the meaning of $c$ might change from line to line
and even in the same chain of inequalities.

\subsection{Statement of the problem and results}
As far as the data of our problem are concerned, let  $\ell$ and  $k_0 >0$ be two real constants.
We also consider the \pier{data} $F$, $\vartheta^* $, $\vartheta_0$ and $\chi_0$ such that
	\begin{gather} \label{R}
F \in  H^1(0, T; H) \cap L^1(0, T; L^{\infty}(\Omega)),
	\\ \label{kappa}
	\vartheta^* \in W, \pier{\quad \vartheta^*>0 \ \textrm{ in $\Omega$},} 
	\\ \label{teta0} 
\vartheta_0 \in V,  \quad \vartheta_0 >0 \ \textrm{ a.e. in $\Omega$}, \quad \ln \vartheta_0 \in H,
	\\ \label{chi0} 
\chi_0 \in W.
	\end{gather}
Moreover, we introduce the functions $\tilde{\beta}$ and $\tilde{\pi}$, satisfying the conditions listed below:
	\begin{gather} \label{beta}
\tilde{\beta}: \mathbb{R} \longrightarrow [0, + \infty] \textrm{ is \pier{lower semicontinuous and} convex with $\tilde{\beta}(0)=0$,} 
		\\ \label{lambda}
	\textrm{$ \tilde{\pi} \in C^1(\mathbb{R})$ and $\pi$ is Lipschitz continuous.}  
	\end{gather}
Since $\tilde{\beta}$ is proper, l.s.c.\ and convex, \revis{its} subdifferential $ \beta:= \partial \tilde{\beta}$ is \revis{a 
well-defined maximal monotone graph.} We denote by 
$D(\beta)$ and $D(\tilde{\beta})$ the effective domains of $\beta$ and $\tilde{\beta}$, respectively. \revis{As} $\tilde{\beta}$ takes on its minimum in $0$, we have that $0 \in \beta(0)$. 
We also assume that
\begin{align} 
\label{chi0b} \chi_0 \in D(\beta) \, \textrm{ a.e. in $ \Omega$, 
and there exists }\, \xi_0 \in H\nonumber \qquad \\
 \textrm{such that} \quad \xi_0 \in \beta (\chi_0) 
 \, \textrm{ a.e. in $ \Omega$},
\end{align} 
whence
\begin{equation} \label{chi0c} 
\tilde{\beta}(\chi_0) \in L^1(\Omega).
\end{equation}
Indeed, thanks to the definition of the subdifferential and to \eqref{beta}, we have that
\begin{equation*} 
0 \leq \int_{\Omega} \tilde{\beta}(\chi_0) \leq (\xi_0, \chi_0) \leq \| \xi_0 \|_H \| \chi_0 \|_H.
\end{equation*}
In the following, the same symbol $\beta$ will be used for the maximal monotone operators induced \colli{by $\beta$} on \pier{$H\equiv L^2(\Omega)$ and $L^2(0,T;H)\equiv L^2(Q) $}.

In our problem\colli{, the} maximal monotone operator
\begin{equation*} \label{A1}
A: H \longrightarrow  H
\end{equation*}
also appears. We assume that 
\begin{align} \label{A1-bis}
A \, \hbox{ is the subdifferential of a convex and l.s.c. function } \, \Phi: H \longrightarrow  \mathbb{R} \nonumber \\
 \quad \textrm{which takes its minimum in $0$ and has at most a quadratic growth}.
\end{align}
\pier{These properties are related to our assumptions on 
$A = \partial \Phi$, which read} 
\begin{equation} 
\pier{0 \in A(0),\quad \exists \, C_A >0 \, \hbox{ such that }\,
\label{stimaA}
\| y \|_H \leq C_A (1 + \| x \|_H )  \quad \forall \,  x \in H, \ \,  \forall \,  y \in A x.}
\end{equation}
In the following, the same symbol $A$ will be used for the maximal monotone \colli{operator} induced on $L^2(0,T;H)$.

\paragraph{Examples of operators $A$.}
Let us consider the operator
\begin{equation*} 
\sign \ : \mathbb{R} \longrightarrow 2^{\mathbb{R}} , \quad 
\label{segno1111111111111}
\sign(r) =  
\left\{ \begin{array}{ll}
\frac{r}{| r |}  								& \textrm{if $r \neq 0$}, \\[0.2cm]
\textrm{$[ - 1, 1 ]$}  										& \textrm{if $ r =0  $,}
\end{array}
\right.
\end{equation*}
and its nonlocal counterpart in $H$, that is, 
\begin{equation*} \label{Sign}
\Sign: H \longrightarrow 2^H , \quad 
\Sign(v) = 
\left\{ \begin{array}{ll}
\frac{v}{\| v \|_H} 								& \textrm{if $v \neq 0$}, \\[0.2cm]
B_1(0)													& \textrm{if $ v =0  $,}
\end{array}
\right.
\end{equation*}
where $B_1(0)$ denotes the closed unit ball of $H$. \pier{It is straightforward to check that 
$\Sign$ satisfies \eqref{A1-bis}--\eqref{stimaA} and turns out to be 
the subdifferential of the norm function $v \mapsto \| v \|_H  $. Concerning the graph
$\sign$, it is well known that it induces a maximal monotone operator in $H$ which is the 
the subdifferential of the convex function $v \mapsto \int_\Omega |v|  $.}

\paragraph{Main result.}
Our aim is to find a \pier{quadruplet $(\vartheta,\chi,\zeta,\xi)$} satisfying the regularity conditions 
	\begin{gather} \label{reg1}
\vartheta \in L^2(0, T; V),
	\\ \label{reg2}
\textrm{$\vartheta >0$ \ a.e. in $Q$ \, and} \quad  \ln{\vartheta}  \in H^1(0, T; V') \cap L^{\infty}(0, T;H),
	\\ \label{reg3}
\chi  \in H^1(0, T; H)\cap C^0([0,T];V) \cap L^2(0, T;W),
	\\ \label{reg4}
\zeta \in  L^2(0, T;H), \quad \quad \xi  \in  L^2(0, T;H),
	\end{gather}
and solving \colli{the} Problem $(P)$ defined by
\begin{gather} \label{iniziale1}
\partial_t( \ln{\vartheta} (t) + \ell \chi(t)) - k_0 \Delta \vartheta(t) + \zeta(t)= F(t) 
\quad \textrm{in $V'$, for a.e. $t \in (0,T)$},
\\ \label{iniziale2}
\partial_t \chi - \Delta \chi + \xi + \pi(\chi) = \ell\vartheta   \quad \textrm{a.e. in $Q$,}
\\ \label{iniziale3}
\zeta(t) \in A(\vartheta(t) - \vartheta^*) \quad \textrm{for a.e. $t \in (0,T)$,}  
\\ \label{iniziale4}
\xi \in \beta(\chi) \quad  \textrm{a.e. in $Q$,}
\\ \label{iniziale5}
\partial_{\nu} \vartheta =0,  \quad \partial_{\nu}\chi = 0 \quad \textrm{in the sense of traces on $\Sigma$},
\\ \label{iniziale6}
\ln{\vartheta}(0) = \ln{\vartheta_0}, \quad \quad \chi(0) = \chi_0 \quad \textrm{a.e. in $\Omega$}.
\end{gather}
Here, we pointed out the boundary conditions \eqref{iniziale5} although they 
are already contained in the specified meaning of $- \Delta$ (cf.~\eqref{pier4}). 
By the way, a variational formulation of \eqref{iniziale1} reads
\begin{align} \label{rifff}
\langle \partial_t( \ln{\vartheta} (t) + \ell \chi(t)) + \zeta(t), v \rangle
+ k_0 \int_{\Omega} \nabla \vartheta (t) \cdot \nabla v 
= \int_{\Omega} F(t) v \nonumber \\
\hbox{for all }\, v \in V, \textrm{ for a.e. $t \in (0,T)$}.
\end{align}
About the initial conditions in \eqref{iniziale6}, note that from \eqref{reg2} it follows that $\ln{\vartheta}$ is at least weakly continuous from $[0,T]$ to $H$. 

The following result is concerned with the existence of solutions to Problem (P).

\begin{theorem} \label{Teo-esistenza}
Assume \eqref{R}--\eqref{stimaA}. 
Then the Problem $(P)$ stated by \eqref{iniziale1}--\eqref{iniziale6}
has at least a solution $(\vartheta,\chi,\zeta,\xi)$ satisfying \eqref{reg1}--\eqref{reg4}
and the regularity properties
\begin{equation} \label{reg1th}
\vartheta \in L^{\infty}(0, T; V), \quad \quad \zeta \in  L^{\infty}(0, T;H),
\end{equation}
\begin{equation} \label{reg2th}
\chi  \in W^{1, \infty}(0, T; H) \cap H^1(0, T;V) \cap L^{\infty}(0, T;W), \quad \quad \xi  \in   L^{\infty}(0, T;H).
\end{equation}
\end{theorem}
The proof of Theorem \ref{Teo-esistenza} will be given in the subsequent three sections. 

\setcounter{equation}{0}
\section{The approximating problem $(P_{\tau})$}

In order to prove the existence theorem, first we introduce a backward finite differences scheme.
Assume that $N$ is a positive integer and let $Z$ be any normed space.
By fixing the time step 
$$\tau = T/N , \quad \colli{N}\in \mathbb{N},$$ 
we introduce the interpolation maps from $Z^{N+1}$
into either $L^{\infty}(0,T; Z)$ or $W^{1,\infty}(0,T; Z)$.
For $(z^0, z^1, \colli{\ldots\,} , z^N) \in Z^{N+1}$, we define
the piecewise constant functions $\overline{z}_{\tau}$
and the piecewise linear functions $\widehat{z}_{\tau}$, respectively:
\begin{align} \label{regtratti1}
&\overline{z}_{\tau}  \in L^{\infty}(0,T; Z), \quad \overline{z}((i+s)\tau)= z^{i+1}, 
\nonumber\\
&\widehat{z}_{\tau} \in W^{1,\infty}(0,T; Z), \quad \widehat{z}((i+s)\tau)= z^i + s(z^{i+1} - z^i), 
\nonumber\\ 
&\hskip5cm \hbox{if $0<s<1$ and $i=0, \colli{\ldots\,}, N-1$.}
\end{align}
By a direct computation, it is straightforward to prove that
\begin{gather} \label{regtratti5}
\| \overline{z}_{\tau} - \widehat{z}_{\tau} \|_{L^{\infty}(0,T; Z)} = \max_{i=0, \colli{\ldots\,}, N-1} \| z_{i+1} - z_i \|_Z 
= \tau \| \partial_t \widehat{z}_{\tau} \|_{L^{\infty}(0,T; Z)},
\\ \label{regtratti6}
\| \overline{z}_{\tau} - \widehat{z}_{\tau} \|^2_{L^2(0,T; Z)} = 
\frac{\tau}{3} \sum_{i=0}^{N-1} \| z_{i+1} - z_i \|^2_Z 
= \frac{\tau^2}{3} \| \partial_t \widehat{z}_{\tau} \|^2_{L^2(0,T; Z)},
\end{gather}
\vskip-.75cm
\begin{align}
\nonumber
\| \overline{z}_{\tau} - \widehat{z}_{\tau} \|^2_{L^{\infty}(0,T; Z)} 
&= \max_{i=0, \colli{\ldots\,}, N-1} \| z_{i+1} - z_i \|^2_Z \\
&\leq \sum_{i=0}^{N-1}  \tau^2 \bigg\|  \frac{z_{i+1} - z_i}{\tau} \bigg\|^2_Z
\leq \tau \| \partial_t \widehat{z}_{\tau} \|^2_{L^2(0,T; Z)}. \label{regtratti7}
\end{align}
Then, we consider the approximating problem $(P_{\tau})$. We set 
\begin{equation} \label{Fapp-i}
F^i:= \frac{1}{\tau} \int_{(i-1) \tau}^{i \tau} F(s) \ ds,   \quad \quad  \textrm{for $i=1, \colli{\ldots\,}, N$,}
\end{equation}
and we look for two vectors 
$(\vartheta^0, \vartheta^1, \colli{\ldots\,} , \vartheta^N) \in V^{N+1}$,
$(\chi^0, \chi^1, \colli{\ldots\,} , \chi^N) \in W^{N+1}$ satisfying, for $i=1, \colli{\ldots\,}, N$, the system
\begin{gather} 
\label{D-ini-0}
\vartheta^i > 0 \quad \hbox{a.e. in $\Omega$}, \quad \ln{\vartheta^i} \in H , \quad 
\exists \ \zeta^i ,\, \xi^i \in H \quad \hbox{such that} 
\\[0.1cm]
\tau^{1/2} \vartheta^i +\ln{\vartheta^i} + \ell \chi^i + \tau \zeta^i -\tau k_0 \Delta \vartheta^i 
= \tau F^i + \tau^{1/2} \vartheta^{i-1} + \ln{\vartheta^{i-1}} + \ell \chi^{i-1} 
\nonumber 
\\ 
\label{D-ini-1} 
\hskip8cm\textrm{ a.e. in $ \Omega$,}
\\[0.1cm]
\label{D-ini-2}
\chi^{i} - \tau \Delta \chi^{i} + \tau \xi^i + \tau \pi(\chi^i) = \chi^{i-1} + \tau \ell \vartheta^i   
\quad  \textrm{ a.e. in $ \Omega$,}
\\ \label{D-ini-3}
\zeta^i \in A(\vartheta^i - \vartheta^*)   ,
\\ \label{D-ini-4}
\xi^i \in \beta(\chi^i) \quad \ \textrm{a.e. in $ \Omega$,}
\\ \label{D-ini-5}
\partial_{\nu} \vartheta^i = \partial_{\nu}\chi^i = 0 \quad  \textrm{ a.e. 
on $\Gamma$},
\\ \label{D-ini-6}
\vartheta^0 = \vartheta_0, \quad  \chi^0 = \chi_0 \quad  \textrm{a.e. in $ \Omega$.}
\end{gather}
\pier{In view of \eqref{R}--\eqref{chi0}, we infer that for $i=1$ the right-hand side of  
\eqref{D-ini-1} is an element of $H$, and \colli{for any given $\chi^1 $ (present in the left-hand side) we have 
to find the corresponding $\vartheta^1$, along with 
$\xi^1 $,} fulfilling \eqref{D-ini-0}--\eqref{D-ini-1} and \eqref{D-ini-3}; in case we succeed, from  a comparison in \eqref{D-ini-1} it will turn out that 
$\vartheta^1 \in W$. 
Then, we insert $\vartheta^1$\colli{, depending on $\chi^1$,} in the right-hand side of  \eqref{D-ini-2} and we seek \colli{somehow a fixed point $\chi^1$ , together with $\xi^1 \in H$,} satisfying \eqref{D-ini-2} and \eqref{D-ini-4}. Once we recover \colli{$\chi^1$ and the related $\vartheta^1$}, we can start again our procedure, and so on. Then, it is important to show that, for a fixed $i $ and known data $F^i, \,  \vartheta^{i-1} ,\, \ln{\vartheta^{i-1}} , \, \chi^{i-1}$ we are able to find a pair $(\vartheta^{i} ,\, \chi^{i}) $ solving \eqref{D-ini-0}--\eqref{D-ini-5}.}

\begin{theorem}\label{teotau}
\pier{There exists some fixed value $\tau_1 \leq \min\{1,T\}$, depending only on the data, 
such that for any time step $0<\tau<\tau_1 $ the approximating problem $(P_{\tau})$
stated by \eqref{D-ini-0}--\eqref{D-ini-6} has a unique solution}
$$
\pier{(\vartheta^0, \vartheta^1, \colli{\ldots\,} , \vartheta^N) \in V \times W^N , \quad \ 
(\chi^0, \chi^1, \colli{\ldots\,} , \chi^N) \in W^{N+1}.}
$$
\end{theorem}
\pier{Let us now rewrite the discrete equation
\eqref{D-ini-1}--\eqref{D-ini-6} by using the piecewise constant and piecewise linear functions defined in \eqref{regtratti1}, 
with obvious notation, and obtain that
\begin{gather} \label{appr1}
\tau^{1/2} \partial_t \widehat{\vartheta}_{\tau} 
+ \partial_t  \widehat{\ln{\vartheta}}_{\tau} 
+ \ell \partial_t  \widehat{\chi}_{\tau}
+ \overline{\zeta}_{\tau}
- k_0 \Delta \overline{\vartheta}_{\tau} = \overline{F}_{\tau} \quad \textrm{a.e. in $Q$,}
\\ \label{appr2}
\partial_t \widehat{\chi}_{\tau} 
- \Delta \overline{\chi}_{\tau} 
+ \overline{\xi}_{\tau} 
+ \pi(\overline{\chi}_{\tau}) =  \ell \overline{\vartheta}_\tau \quad \textrm{a.e. in $Q$,}
\\ \label{appr3}
\overline{\zeta}_{\tau} (t) \in A(\overline{\vartheta}_{\tau} (t) - \vartheta^*)  
 \quad \textrm{for a.e. $t \in (0,T)$,}
\\ \label{appr4}
\overline{\xi}_{\tau} \in \beta(\overline{\chi}_{\tau}) \quad \textrm{a.e. in $Q$,}
\\ \label{appr5}
\partial_{\nu} \overline{\vartheta}_{\tau} = \partial_{\nu} \overline{\chi}_{\tau} = 0 \quad \textrm{a.e. on $\Sigma$},
\\ \label{appr6}
\widehat{\vartheta}_{\tau} \, (0) = \vartheta_0, \quad \widehat{\chi}_{\tau} \,(0) = \chi_0 \quad \textrm{a.e. in $ \Omega$.}
\end{gather}}

\subsection{The auxiliary approximating problem $(AP_{\varepsilon})$}
In this subsection we \colli{introduce} the auxiliary 
approximating problem $(AP_{\varepsilon})$
obtained by considering the approximating problem $(P_{\tau})$ at each step $i=1,\ldots, N$ and replacing the monotone operators appearing in \eqref{D-ini-0}--\eqref{D-ini-6}
with their Yosida regularizations. \pier{About general properties of maximal monotone operators and subdifferentials of convex functiions, we refer the reader to \cite{Barbu, Brezis}.}

\paragraph{Yosida regularization of $ \ln $.} 
We introduce the Yosida regularization of $\ln$. For $\varepsilon>0$ we set
\begin{equation} \label{richiamolog1}
{\ln}_{\varepsilon}: \mathbb{R} \longrightarrow \mathbb{R},
\quad \ 
{\ln}_{\varepsilon} : = \frac{I - (I + \varepsilon {\ln})^{-1}}{\varepsilon}. 
\end{equation}
where $I$ denotes the identity. We \revis{point out} that ${\ln}_{\varepsilon}$ is monotone, 
Lipschitz continuous (with Lipschitz constant $1/\varepsilon$) and satisfies the 
following properties: denoting by $L_{\varepsilon}= (I + \varepsilon {\ln})^{-1}$ the 
resolvent operator, we have that
	\begin{gather*}
	{\ln}_{\varepsilon}(x) \in  {\ln}(L_{\varepsilon}x) \quad \hbox{ for all } \, x\in \mathbb{R},\\
	|{\ln}_{\varepsilon}(x)| \leq |{\ln}(x)|,
\quad   	\lim_{\varepsilon \searrow 0} { \ln }_{\varepsilon}(x) =  { \ln }(x) \quad 
\hbox{ for all } \, x > 0.	
\end{gather*}
\pier{We also introduce the nonnegative and convex functions 
\begin{equation}\label{defLam}
\Lambda (x)= \int_1^x \ln r \, dr , \quad \Lambda_{\varepsilon} (y)= \int_1^y \ln_{\varepsilon} r \, dr \quad \hbox{ for all } \, x > 0 \hbox{ and } y\in \mathbb{R}.
\end {equation}
Note that the  graph $x \mapsto \ln x $ is nothing but the subdifferential of the convex function $
\Lambda$ extended by lower semicontinuity in $0$ and with value $+\infty $ for $x<0$.
On the other hand, $\Lambda_{\varepsilon}$ coincides with the Moreau--Yosida regularization of $\Lambda$ and, in particular, we have that
\begin{equation} \label{propLam}
	0\leq \Lambda_{\varepsilon}(x) \leq \Lambda(x) \quad \textrm{for every $x >0$.}
	\end{equation}
}

\paragraph{Yosida regularization of $A$.} 
We introduce the Yosida regularization of $A$. For $\varepsilon>0$ we define
\begin{equation} \label{richiamoA1}
A_{\varepsilon}: H \longrightarrow  H, \quad \quad
A_{\varepsilon} = \frac{I - (I + \varepsilon A)^{-1}}{\varepsilon}.
\end{equation} 
Note that $A_{\varepsilon}$ is Lipschitz-continuous (with Lipschitz
constant $1/\varepsilon$) and maximal monotone in $H$. Moreover, 
$A$ satisfies the following properties:
denoting by $J_{\varepsilon}= (I + \varepsilon A)^{-1}$ the resolvent operator, for all $\delta > 0$
and for all $x \in H$, we have that
	\begin{gather} \label{AJe}
	A_{\varepsilon}x \in A(J_{\varepsilon}x),
	\\[0.2cm]
	\label{pier5}	
	\| A_{\varepsilon}x \|_H \leq \| \colli{A^\circ} x \|_H,
	\quad 
	\lim_{\varepsilon \searrow 0} \| A_{\varepsilon}x - \colli{A^\circ} x \|_H =0, 
	\end{gather}
where $\colli{A^\circ} x$ is the element	of the range of $A$ having minimal norm. 
\pier{Let us point out a key property of $A_{\varepsilon}$, 
which is a consequence of \eqref{stimaA}: indeed, there holds
\begin{equation} \label{Ae}
\| A_{\varepsilon} x \|_H \leq C_A (1 + \| x \|_H ) \quad 
\textrm{ for all} \ \, x \in H.
\end{equation}
Notice that $0 \in A(0)$ and $0 \in I(0)$: consequently, for every $\varepsilon >0$ 
we infer that $ J_\varepsilon (0)=0$.}
Moreover, since $A$ is  maximal monotone,  $J_\varepsilon$ is a contraction. Then, from \eqref{stimaA} and \eqref{AJe} it follows that
\begin{align*} \nonumber
\| A_{\varepsilon}x \|_H  \leq C_A (\| J_{\varepsilon}x\|_H  + 1)  \leq  C_A (\| J_{\varepsilon}x - J_{\varepsilon}0 \|_H  + 1) \leq  C_A(\| x \|_H + 1) \\
\quad \hbox{ for every $x\in H$}.
\end{align*}

\paragraph{Yosida regularization of $\beta$.} 
We introduce the Yosida regularization of $\beta$. For $\varepsilon>0$~\revis{let}
\begin{equation} \label{richiamobeta1}
\beta_{\varepsilon}: \mathbb{R} \longrightarrow \mathbb{R},
\quad \quad
\beta_{\varepsilon} = \frac{I - (I + \varepsilon \beta)^{-1}}{\varepsilon}. 
\end{equation}
\pier{We remark that $\beta_{\varepsilon}$ is Lipschitz continuous (with Lipschitz constant $1/\varepsilon$) and satisfies the following properties: denoting by $R_{\varepsilon}= (I + \varepsilon \beta)^{-1}$ the resolvent operator, we have that
	\begin{gather*}
	\beta_{\varepsilon}(x) \in \beta(R_{\varepsilon}x) \quad \hbox{ for all } \, x\in \mathbb{R},\\
	|\beta_{\varepsilon}(x)| \leq |\colli{\beta^\circ}(x)|, \quad
	\lim_{\varepsilon \searrow 0} \beta_{\varepsilon}(x) =  \colli{\beta^\circ}(x) \quad 
\hbox{ for all } \, x \in D(\beta),	
\end{gather*}
where $\colli{\beta^\circ}(x)$ is the element	of the range of $\beta (x)$ having minimal modulus.
We also introduce the Moreau--Yosida} regularization of $\tilde{\beta}$. For $\varepsilon>0$ 
and $x\in \mathbb{R}$ we \revis{set}
\begin{equation*}
\tilde{\beta}_{\varepsilon}: \mathbb{R} \longrightarrow [0, +\infty],
\quad \quad \pier{\tilde{\beta}_{\varepsilon}(x) := \min_{y \in \mathbb{R}} 
\left\{\tilde{\beta}(y) + \frac{1}{2 \varepsilon} |x - y|^2\right\} }
\end{equation*}
and recall that
	\begin{equation} \nonumber 
	\tilde{\beta}_{\varepsilon}(x) \leq \tilde{\beta}(x) \quad \textrm{for every $x  \in \mathbb{R}$.}
	\end{equation}
We also observe that $\beta_{\varepsilon}$ is \pier{the derivative of $\tilde{\beta}_{\varepsilon}$. Then, for every $x_1, x_2  \in \mathbb{R}$} we have that
\begin{equation} \nonumber 
\tilde{\beta}_{\varepsilon}(x_2) = \tilde{\beta}_{\varepsilon}(x_1) + \int_{x_1}^{x_2} \beta_{\varepsilon}(s) \ ds.
\end{equation}

\paragraph{Definition of the auxiliary approximating problem $(AP_{\varepsilon})$.}
We fix $\tau $ and specify an auxiliary approximating problem $(AP_{\varepsilon})$, which
is obtained by considering \eqref{D-ini-0}--\eqref{D-ini-5} for a fixed $i$ and 
introducing the regularized operators defined above. We set
\begin{equation} \label{gh}
g: = \tau F^i + \tau^{1/2} \vartheta^{i-1} + \ln{\vartheta^{i-1}} + \ell \chi^{i-1},
\quad \quad h: = \chi^{i-1}, 
\end{equation}
\pier{and note that both $g$ and $h$ are prescribed elements of $H$ (cf. \eqref{Fapp-i}, \eqref{R}, \eqref{teta0}, \eqref{chi0} and \eqref{D-ini-0}).}
We look for a pair $(\Theta_{\varepsilon},X_{\varepsilon})$ such that
\begin{gather} \label{E-D-ini-1}
\tau^{1/2} \Theta_{\varepsilon} + {\ln}_{\varepsilon}{\Theta_{\varepsilon}} + \tau A_{\varepsilon} (\Theta_{\varepsilon} - \vartheta^*) 
-\tau k_0 \Delta \Theta_{\varepsilon}
= - \ell X_{\varepsilon} + g \quad  \textrm{ a.e. in $\Omega$,}
\\ \label{E-D-ini-2}
X_{\varepsilon} - \tau \Delta X_{\varepsilon} + \tau \beta_{\varepsilon}(X_{\varepsilon}) 
+ \tau \pi(X_{\varepsilon}) = h + \tau \ell \Theta_{\varepsilon} 
\quad \textrm{ a.e. in $ \Omega$,}
\end{gather}
where $\ln_{\varepsilon}$,  $A_{\varepsilon}$ and $\beta_{\varepsilon}$ are the Yosida regularization of 
$\ln$, $A$ and $\beta$ defined by \eqref{richiamolog1}, \eqref{richiamoA1} and \eqref{richiamobeta1}, respectively. \pier{Here, according to the extended meaning of $ -\Delta $ (see~\eqref{pier4}), we omit the specification of the boundary conditions as with \eqref{D-ini-5}.}

\begin{theorem}\label{teotau2}
\pier{Let $g,\, h \in H$. Then there exists some fixed value $\tau_2 \leq \min\{1,T\}$, depending only on the data, such that for every time step $\tau \in (0,\tau_2)$ and for all $\varepsilon \in (0,1] $ the auxiliary approximating problem $(AP_{\varepsilon})$ 
stated by \eqref{E-D-ini-1}--\eqref{E-D-ini-2} has a unique solution~$(\Theta_{\varepsilon},X_{\varepsilon})$.}
\end{theorem}

\subsection{Existence of a solution for $(AP_{\varepsilon})$}
In order to prove the existence of the solution for the auxiliary 
approximating problem $(AP_{\varepsilon})$
we intend to apply \pier{\cite[Corollary 1.3, p.~48]{Barbu}.}
\pier{To this aim, we point out that, for $\tau$ small enough, the two operators 
\begin{align} \label{oper1}
[\tau^{1/2} I + {\ln}_{\varepsilon} + \tau A_{\varepsilon} (\, \cdot \,  - \vartheta^*) -\tau k_0 \Delta ]  \quad  \hbox{appearing in \eqref{E-D-ini-1},}
\\ \label{oper2}
[I + \tau \beta_{\varepsilon} + \tau \pi - \tau \Delta]  \quad
 \hbox{appearing in \eqref{E-D-ini-2}},
\end{align}
both with domain $W$ and range $H$, are maximal monotone and coercive. 
Indeed, they are the sum of a monotone, Lipschitz continuous and coercive operator:
$$ \tau^{1/2} I + {\ln}_{\varepsilon} + \tau A_{\varepsilon} (\, \cdot \,  - \vartheta^*)
\quad \hbox{in \eqref{oper1}, \ and} \quad I + \tau \beta_{\varepsilon} + \tau \pi \quad \hbox{in \eqref{oper2}}, $$
and of a maximal monotone operator that is $- \Delta$ with a positive coefficient in front. We now check our first claim. Letting $v_1$, $v_2 \in H$, we have that
\begin{align} \nonumber
&\left( (\tau^{1/2} I + {\ln}_{\varepsilon} + \tau A_{\varepsilon} (\, \cdot \,- \vartheta^*))(v_1) 
- (\tau^{1/2} I + {\ln}_{\varepsilon} + \tau A_{\varepsilon} (\, \cdot \,- \vartheta^*) )(v_2), v_1 - v_2 \right)
\\ 
&{}\geq \tau^{1/2} \| v_1 - v_2 \|^2_H + \big({\ln}_{\varepsilon}(v_1) - {\ln}_{\varepsilon}(v_2), v_1 - v_2 \big) \nonumber
\\ &\quad
+ \tau \big(A_{\varepsilon}(v_1 - \vartheta^* ) - A_{\varepsilon}(v_2 - \vartheta^*), (v_1 - \vartheta^*) - (v_2 -  \vartheta^*) \big). \nonumber
\end{align}
Due to the monotonicity of ${\ln}_{\varepsilon}$ and $A_{\varepsilon}$, we have that the 
last two terms on the right-hand side are nonnegative, so that  
\begin{align} \nonumber
&\left( (\tau^{1/2} I + {\ln}_{\varepsilon} + \tau A_{\varepsilon} (\, \cdot \,- \vartheta^*))(v_1) 
- (\tau^{1/2} I + {\ln}_{\varepsilon} + \tau A_{\varepsilon} (\, \cdot \,- \vartheta^*) )(v_2), v_1 - v_2 \right) \\
&\geq \tau^{1/2} \| v_1 - v_2 \|^2_H, \label{monop21}
\end{align}
i.e., the operator $\tau^{1/2} I + {\ln}_{\varepsilon} + \tau A_{\varepsilon} (\, \cdot \,  - \vartheta^*)$
is strongly monotone, hence coercive, in $H$. 
Next, for all $v_1$, $v_2 \in H$ we have that
\begin{align} \nonumber
&\left( (I + \tau \beta_{\varepsilon} + \tau \pi)(v_1) - (I + \tau \beta_{\varepsilon} + \tau \pi)(v_2), v_1 - v_2 \right)
\\ \label{mon1}
&\geq \| v_1 - v_2 \|^2_H + \tau (\beta_{\varepsilon}(v_1) - \beta_{\varepsilon}(v_2), v_1 - v_2 )
- C_{\pi} \tau \| v_1 - v_2 \|^2_H.
\end{align}
where $C_{\pi}$ denotes a Lipschitz constant for $\pi$.  
Since $\beta_{\varepsilon}$ is monotone, it turns out that
\begin{equation*}
\tau (\beta_{\varepsilon}(v_1) - \beta_{\varepsilon}(v_2), v_1 - v_2 ) \geq 0
\end{equation*}
and, choosing $ \tau_2 \leq 1 / 2C_{\pi}$, from \eqref{mon1} we infer that
\begin{equation} \label{maxmonop2}
( (I + \tau \beta_{\varepsilon} + \tau \pi)(v_1) - (I + \tau \beta_{\varepsilon} + \tau \pi)(v_1), v_1 - v_2 )
\geq \frac{1}{2} \| v_1 - v_2 \|^2_H,
\end{equation}
whence the operator $I + \tau \beta_{\varepsilon} + \tau \pi$ is strongly monotone and coercive in $H$, for every $\tau \leq \tau_2$.} 

\pier{Now, in order to prove Theorem~\ref{teotau2},} we divide the proof into two steps.
In the first step, we fix $\overline{\Theta}_{\varepsilon} \in H$ in place of 
$\Theta_{\varepsilon} $ on the right-hand side of \eqref{E-D-ini-2}
and find a solution $X_{\varepsilon}$ for \eqref{E-D-ini-2}.
In the second step, we insert on the right-hand side of \eqref{E-D-ini-1} the element $X_{\varepsilon}$ 
obtained in the first step and find a solution $\Theta_{\varepsilon}$ to \eqref{E-D-ini-1}.
Now, let $\overline{\Theta}_{1, \varepsilon}$ and $\overline{\Theta}_{2, \varepsilon}$ be two different input data. 
We denote by $X_{1, \varepsilon}$, $X_{2, \varepsilon}$ the corresponding solutions for \eqref{E-D-ini-2} obtained in the first step 
and by $\Theta_{1, \varepsilon}$, $\Theta_{2, \varepsilon}$ 
the related solution of \eqref{E-D-ini-1} found in the second step.

Hence, taking the difference between the two equations \eqref{E-D-ini-2} written for 
$\overline{\Theta}_{1, \varepsilon}$ and $\overline{\Theta}_{2, \varepsilon}$
and testing the \pier{result} by $(X_{1, \varepsilon} - X_{2, \varepsilon})$, we have that
\begin{align} \nonumber
&\big( (I + \tau \beta_{\varepsilon} + \tau \pi)(X_{1, \varepsilon}) 
- (I + \tau \beta_{\varepsilon} + \tau \pi)(X_{2, \varepsilon} ), 
X_{1, \varepsilon} - X_{2, \varepsilon}  \big) 
\\ \label{intermezz1}
&+ \tau \int_{\Omega} | \nabla (X_{1, \varepsilon} - X_{2, \varepsilon})|^2
\leq \tau \ell \big(\overline{\Theta}_{1, \varepsilon} - \overline{\Theta}_{2, \varepsilon}, X_{1, \varepsilon} - X_{2, \varepsilon} \big).
\end{align}
Then, applying \eqref{maxmonop2} and \eqref{dis2} to the first term on the left-hand side of \eqref{intermezz1}
and to the right-hand side of \eqref{intermezz1}, respectively, we infer that
\begin{equation} \nonumber
\frac{1}{2} \| X_{1,\varepsilon} - X_{2, \varepsilon} \|^2_H 
+ \tau \int_{\Omega} | \nabla (X_{1, \varepsilon} - X_{2, \varepsilon})|^2
\leq \frac{1}{4} \| X_{1, \varepsilon} - X_{2, \varepsilon} \|^2_H 
+ \tau^2 \ell^2 \| \overline{\Theta}_{1, \varepsilon} - \overline{\Theta}_{2, \varepsilon} \|^2_H,
\end{equation} 
whence
\begin{equation} \label{ultimastep1}
\| X_{1, \varepsilon} - X_{2, \varepsilon} \|^2_H 
\leq 4 \tau^2 \ell^2 \| \overline{\Theta}_{1, \varepsilon} - \overline{\Theta}_{2, \varepsilon} \|^2_H.
\end{equation} 
Now, we take the difference between the corresponding equations \eqref{E-D-ini-1} 
written for the solutions $X_{1,\varepsilon}$, $X_{2,\varepsilon}$ obtained in the first step 
and test by $(\Theta_{1, \varepsilon} - \Theta_{2, \varepsilon})$. We obtain that
\begin{align} \nonumber
&\left( (\tau^{1/2} I + {\ln}_{\varepsilon} + \tau A_{\varepsilon} (\, \cdot \,- \vartheta^*))(\Theta_{1, \varepsilon}) 
- (\tau^{1/2} I + {\ln}_{\varepsilon} + \tau A_{\varepsilon} (\, \cdot \,- \vartheta^*) )(\Theta_{2, \varepsilon} ), \Theta_{1, \varepsilon} - \Theta_{2, \varepsilon}  \right)
\\
&+ \tau k_0 \int_{\Omega} | \nabla (\Theta_{1,\varepsilon} - \Theta_{2,\varepsilon})|^2
\label{ultimastep2}
\leq \frac{\ell^2}{2 \tau^{1/2}} \| X_{1, \varepsilon} - X_{2, \varepsilon} \|^2_H  + \frac{\tau^{1/2}}{2} 
\| \Theta_{1, \varepsilon} - \Theta_{2, \varepsilon} \|^2_H.
\end{align}
\pier{By recalling \eqref{monop21} and using it in the left-hand side of \eqref{ultimastep2}
we infer that 
$$ \tau^{1/2} 
\| \Theta_{1, \varepsilon} - \Theta_{2, \varepsilon} \|^2_H \leq
 \frac{\ell^2}{\tau^{1/2}} \| X_{1, \varepsilon} - X_{2, \varepsilon} \|^2_H  $$
Then, by combining this inequality with \eqref{ultimastep1}, we deduce that
\begin{equation} \label{contrazione1}
 \| \Theta_{1, \varepsilon} - \Theta_{2, \varepsilon} \|^2_H
\leq 4 \tau \ell^4 \| \overline{\Theta}_{1, \varepsilon} - \overline{\Theta}_{2, \varepsilon} \|^2_H,
\end{equation}
whence we obtain a contraction mapping for every $\tau \leq \tau_2$, 
provided that $\tau_2 \leq 1/(8 \ell^4)$. Finally, by applying the Banach 
fixed point theorem, we conclude that there exists a unique solution $(\Theta_{\varepsilon}, X_{\varepsilon})$ to the auxiliary problem $(AP_{\varepsilon})$.}

\subsection{A priori estimates on $AP_{\varepsilon}$}
In this subsection we \pier{derive} a series of a priori estimates, independent of $\varepsilon$,
inferred from the equations \eqref{E-D-ini-1}--\eqref{E-D-ini-2}
of the auxiliary approximating problem $(AP_{\varepsilon})$.

\paragraph{First a priori estimate.}
We test \eqref{E-D-ini-1} \pier{by $\tau (\Theta_{\varepsilon} - \vartheta^*)$ and \eqref{E-D-ini-2} by $X_{\varepsilon}$, then we sum up. By exploiting the 
cancellation of the suitable corresponding terms and recalling the definition \eqref{defLam} of $\Lambda_{\varepsilon}$,} we obtain that
\begin{align} \nonumber
&\tau^{3/2} \| \Theta_{\varepsilon} - \vartheta^* \|^2_H 
+ \tau \Lambda_{\varepsilon} (\Theta_{\varepsilon})
\colli{{}+ \tau^2 ( A_{\varepsilon} (\Theta_{\varepsilon} - \vartheta^*), \Theta_{\varepsilon} - \vartheta^* ) + \tau^2 k_0 \int_{\Omega} |\nabla (\Theta_{\varepsilon} - \vartheta^*)|^2{}}\\ \nonumber
&{}\colli{{}+\big( (I + \tau \beta_{\varepsilon} + \tau \pi)(X_{\varepsilon}) 
- (I + \tau \beta_{\varepsilon} + \tau \pi)(0 ), X_{\varepsilon}\big)
+ \tau \int_{\Omega} |\nabla X_{\varepsilon}|^2{}}
\\ \nonumber
&{}\leq
{}- \tau^{3/2} ( \vartheta^*, \Theta_{\varepsilon} - \vartheta^*)
+ \tau \Lambda_{\varepsilon} (\vartheta^*)
{}- \tau^2 k_0 \int_{\Omega} \nabla \vartheta^* 
\cdot \nabla (\Theta_{\varepsilon} - \vartheta^*)
\\ &\quad \ {}+ \ell \tau (X_{\varepsilon}, \vartheta^*)
+ \tau (g, \Theta_{\varepsilon} - \vartheta^*) 
\colli{{}- \tau (\pi(0), X_{\varepsilon} )}
+  (h, X_{\varepsilon}).
\label{stma-1}
\end{align}
Let us note that \colli{all terms on the left-hand side are nonnegative; in particular, recalling \eqref{maxmonop2}, we have that 
\begin{equation} \label{stma-1au1}
\big( (I + \tau \beta_{\varepsilon} + \tau \pi)( X_{\varepsilon}) - (I + \tau \beta_{\varepsilon} + \tau \pi)(0),  X_{\varepsilon} \big)
\geq \frac{1}{2} \|  X_{\varepsilon} \|^2_H, 
\end{equation}
Due} to \eqref{kappa} and the continuity of the positive function $ \vartheta^*$, \eqref{propLam} helps us in estimating the second term on the right-hand side of \eqref{stma-1}: 
\begin{equation} 
\tau \Lambda_{\varepsilon} (\vartheta^*) \leq \tau \Lambda (\vartheta^*) \leq c \,\tau.
\end{equation}
Since $g,\, h  \in H$ and \eqref{kappa} holds, by applying the Young inequality \eqref{dis2} to the \pier{other terms} on the right-hand side of \eqref{stma-1}, we find that
\begin{gather}
- \tau^{3/2} ( \vartheta^*, \Theta_{\varepsilon} - \vartheta^*) \leq \frac{\tau^{3/2}}{4} \| \Theta_{\varepsilon} - \vartheta^* \|^2_H 
+c \,\tau^{3/2},
\\
- \tau^2 k_0 \int_{\Omega} \nabla \vartheta^* \cdot \nabla (\Theta_{\varepsilon} - \vartheta^*) \leq
\frac{\tau^2 k_0}{2} \int_{\Omega} |\nabla (\Theta_{\varepsilon} - \vartheta^*)|^2 + c\,\tau^2,
\\
\ell \tau (X_{\varepsilon}, \vartheta^*) \leq \colli{\frac{1}{8}} \| X_{\varepsilon} \|^2_H  + c\,\tau^2,
\end{gather}
\begin{gather}
\tau (g, \Theta_{\varepsilon} - \vartheta^*) \leq \frac{\tau^{3/2}}{4} \| \Theta_{\varepsilon} - \vartheta^* \|^2_H+ c\,\tau^{1/2},
\\ \label{stma-1au2}
\colli{{}- \tau (\pi(0), X_{\varepsilon} ) \leq  \frac{1}{8} \| X_{\varepsilon} \|^2_H
+ c \, \tau^2 ,}
\quad \quad  (h, X_{\varepsilon}) \leq \colli{\frac{1}{8}} \| X_{\varepsilon} \|^2_H + c.
\end{gather} 
Then, in view of \eqref{stma-1au1}--\eqref{stma-1au2}, \pier{from \eqref{stma-1} and \eqref{kappa} it is not difficult to} infer that 
\begin{equation} \label{stima-1}
\tau^{3/4} \| \Theta_{\varepsilon} \|_H + \tau \| \nabla \Theta_{\varepsilon} \|_H + 
\| X_{\varepsilon} \|_H + \tau^{1/2} \| \nabla X_{\varepsilon} \|_H \leq c
\end{equation} 
taking into account that \colli{$\tau \leq \tau_2$}.

\paragraph{Second a priori estimate.}
We test \eqref{E-D-ini-2} by $\beta_{\varepsilon}(X_{\varepsilon})$ and obtain that
\begin{align} \nonumber
&\big(X_{\varepsilon} , \beta_{\varepsilon}(X_{\varepsilon}) \big) + \tau \int_{\Omega} \beta'_{\varepsilon}(X_{\varepsilon}) |\nabla X_{\varepsilon}|^2
+ \tau \int_{\Omega} |\beta_{\varepsilon}(X_{\varepsilon})|^2
\\ \label{stma-2}
&\leq {}- \tau  \int_{\Omega} \pi(X_{\varepsilon}) \beta_{\varepsilon}(X_{\varepsilon})
+ \tau \ell \int_{\Omega} \Theta_{\varepsilon} \beta_{\varepsilon}(X_{\varepsilon})
+ \int_{\Omega} h \beta_{\varepsilon}(X_{\varepsilon}).
\end{align}
Thanks to the monotonicity of $\beta_{\varepsilon}$ and to the condition 
$\beta_{\varepsilon}(0)=0$, the terms on the left-hand side are nonnegative. 
As $\pi$ is Lipschitz continuous, by applying the Young inequality \eqref{dis2} 
to every term on the right-hand side of \eqref{stma-2} and using \eqref{stima-1}, \pier{for $0< \tau \leq 1$ we obtain that}
\begin{gather} \label{grill1}
 - \tau  \int_{\Omega} \pi(X_{\varepsilon}) \beta_{\varepsilon}(X_{\varepsilon})  
\leq \frac{\tau}{4} \int_{\Omega} |\beta_{\varepsilon}(X_{\varepsilon})|^2 + c,
\\ \label{grill2}
\tau \ell \int_{\Omega} \Theta_{\varepsilon} \beta_{\varepsilon}(X_{\varepsilon})
\leq \frac{\tau}{4} \int_{\Omega} |\beta_{\varepsilon}(X_{\varepsilon})|^2 + \pier{\frac{c}{\tau^{1/2}}},
\\ \label{grill3}
\int_{\Omega} h \beta_{\varepsilon}(X_{\varepsilon}) \leq \frac{\tau}{4} \int_{\Omega} |\beta_{\varepsilon}(X_{\varepsilon})|^2 + \pier{\frac{c}{\tau}}.
\end{gather}
Then, \pier{owing to \eqref{grill1}--\eqref{grill3}, from \eqref{stma-2} it follows that 
\begin{equation}\label{stima-12}
\tau \|\beta_{\varepsilon}(X_{\varepsilon})\|_H^2 \leq \pier{c\big( 1 + \tau^{-1} \big)},
\quad \hbox{ so that } \quad \tau \|\beta_{\varepsilon}(X_{\varepsilon})\|_H \leq c. 
\end{equation} 
Hence, by comparison in \eqref{E-D-ini-2}, we conclude that $\tau \| \Delta X_{\varepsilon} \|_H \leq c$ and, from \eqref{stima-1} and standard elliptic regularity results,}
\begin{equation} \label{stima-13}
\tau \| X_{\varepsilon} \|_W \leq c.
\end{equation}

\paragraph{Third a priori estimate.}
\pier{Recalling \eqref{Ae}, \eqref{kappa} and \eqref{stima-1}, we immediately deduce that
\begin{equation}
\label{pier6} 
\tau \| A_{\varepsilon} (\Theta_{\varepsilon} - \vartheta^*) \|_H \leq \tau\, C_A (1+ 
 \|\Theta_{\varepsilon}\|_H  + \| \vartheta^* \|_H ) \leq c. 
\end{equation} 
Next,} we test \eqref{E-D-ini-1} by ${\ln}_{\varepsilon}{\Theta_{\varepsilon}}$ and obtain that
\begin{align} \nonumber
\| {\ln}_{\varepsilon}{\Theta_{\varepsilon}}\|^2_H + \tau k_0 \int_{\Omega} {\ln}'_{\varepsilon} (\Theta_{\varepsilon}) 
|\nabla \Theta_{\varepsilon}|^2
\leq  {}- \tau^{1/2} (\Theta_{\varepsilon}, {\ln}_{\varepsilon}{\Theta_{\varepsilon}})
\\ \nonumber
{}- \tau \big( A_{\varepsilon} (\Theta_{\varepsilon} - \vartheta^*), {\ln}_{\varepsilon}{\Theta_{\varepsilon}} \big)
- \ell ( X_{\varepsilon}, {\ln}_{\varepsilon}{\Theta_{\varepsilon}}) + (g, {\ln}_{\varepsilon}{\Theta_{\varepsilon}}). 
\end{align}
Then, by applying the \pier{Cauchy--Schwarz inequality to every term 
on the right-hand side and using \eqref{stima-1} and \eqref{pier6}, we infer} that
\begin{equation*}
\| {\ln}_{\varepsilon}{\Theta_{\varepsilon}}\|_H \leq \pier{{}\tau^{1/2} \| \Theta_{\varepsilon} \|_H + c \leq c \big(\tau^{-1/4} + 1 \big),}
\end{equation*}
whence 
\begin{equation} \label{stima-14}
\tau^{1/4} \| {\ln}_{\varepsilon}{\Theta_{\varepsilon}}\|_H \leq  c.
\end{equation}
Moreover, due to \pier{\eqref{stima-14} and \eqref{stima-1}, by comparison in \eqref{E-D-ini-1} it is straightforward to see that
$\tau^{5/4} \| \Delta \Theta_{\varepsilon} \|_H \leq  c$ and consequently}
\begin{equation} \label{stima-15}
\tau^{5/4} \| \Theta_{\varepsilon}\|_W \leq  c.
\end{equation}

\subsection{Passage to the limit as $\varepsilon \searrow 0$}
In this subsection we pass to the limit as $\varepsilon \searrow 0$ and prove that
the limit of subsequences of solutions $(\Theta_{\varepsilon},X_{\varepsilon})$
for ($AP_{\varepsilon}$) (see \eqref{E-D-ini-1}--\eqref{E-D-ini-2}) yields a solution $(\vartheta^i,\chi^i)$ \pier{to \eqref{D-ini-0}--\eqref{D-ini-4}; then, we can conclude that the problem ($P_{\tau}$) has a solution.}

Since the constants appearing in \eqref{stima-1} and
\eqref{stima-12}--\eqref{stima-15} 
do not depend on $\varepsilon$,
we infer that, at least for a subsequence, there exist some limit functions \pier{$(\vartheta^i, \chi^i, L^i, Z^i, B^i)$} such~that
\begin{align} \label{ld1}
&\Theta_{\varepsilon}       \rightharpoonup  \vartheta^i     \quad \hbox{and}\quad 
X_{\varepsilon}            \rightharpoonup  \chi^i    \quad  \hbox{ in } \,  W, \\ 
&{\ln}_{\varepsilon} (\Theta_{\varepsilon} )     \rightharpoonup  L^i,\quad 
 A_{\varepsilon} (\Theta_{\varepsilon} - \vartheta^*)       \rightharpoonup  Z^i 
 \quad \hbox{and}\quad  \beta_{\varepsilon} (X_{\varepsilon})         \rightharpoonup  B^i   \quad     \hbox{ in } \, H, \label{ld4}
 \end{align} 
as $ \varepsilon \searrow 0$. Thanks to the well-known compact embedding $W\subset V$,  
from \eqref{ld1} we infer~that
\begin{equation} \label{ldf1}
\Theta_{\varepsilon}     \rightarrow  \vartheta^i   \quad \hbox{and}\quad 
X_{\varepsilon}          \rightarrow  \chi^i      \quad  \hbox{ in } \,  V.
\end{equation} 
\pier{Besides}, as $\pi$ is Lipschitz continuous, we have that
$|\pi(X_{\varepsilon}) - \pi(\chi^i)| \leq  C_{\pi} |X_{\varepsilon} - \chi^i|$,
whence,  thanks to \eqref{ldf1}, we obtain that
\begin{equation} \label{cvgpigr}
\pi(X_{\varepsilon}) \rightarrow \pi(\chi^i) \quad \quad \textrm{in $H$,}
\end{equation}
as $ \varepsilon \searrow 0$. 
Now, we pass to the limit on \pier{${\ln}_{\varepsilon}(\Theta_{\varepsilon})$,  
$A_{\varepsilon}(\Theta_{\varepsilon} - \vartheta^*)$ and  $\beta_{\varepsilon}
(X_{\varepsilon})$. In view of a general convergence result involving maximal monotone 
operators (see, e.g., \cite[Proposition~1.1, p.~42]{Barbu}), thanks 
 to the strong convergences in $H$ ensured by \eqref{ldf1} and to the weak convergences in \eqref{ld4},} we conclude that
\begin{equation} \label{convoperatorivari}
 L^i    \in {\ln}(\chi^i),  
 \quad \  Z^i    \in A(\vartheta^i - \vartheta^*), 
 \ \quad B^i    \in \beta(\chi^i). 
 \end{equation}
\pier{In conclusion}, using \eqref{ld1}--\eqref{convoperatorivari} and recalling \eqref{gh}, we can pass to the limit as $\varepsilon \searrow 0$ 
in \eqref{E-D-ini-1}--\eqref{E-D-ini-2} so to
obtain \eqref{D-ini-0}--\eqref{D-ini-4} for the limiting functions $\vartheta^i$ and $\chi^i$.

\subsection{Uniqueness of the solution of $(P_{\tau})$}
In this section we prove that the approximating problem $(P_{\tau})$ stated by
\pier{\eqref{D-ini-0}--\eqref{D-ini-6} has a unique solution. Then, the proof of Theorem~\ref{teotau} will be complete.}

We write problem $(P_{\tau})$ for two solutions $(\vartheta^i_1, \chi^i_1)$, $(\vartheta^i_2, \chi^i_2)$
and set $\vartheta^i := \vartheta^i_1 - \vartheta^i_2$ and $\chi^i:= \chi^i_1 - \chi^i_2$, \colli{$i= 1,\ldots, N.$}
Then, we multiply by $\tau \vartheta^i$ the difference between the corresponding equations \eqref{D-ini-1}
and by $\chi^i$ the difference between the corresponding equations \eqref{D-ini-2}.
Adding the resultant equations, we obtain that
\begin{align} \nonumber
&\tau^{3/2} \| \vartheta^i \|^2_H 
+ \tau \big( \ln{\vartheta^i_1} - \ln{\vartheta^i_2}, \vartheta^i_1 - \vartheta^i_2 \big)
+ \tau^2 \big(\zeta^i_1 - \zeta^i_2, \vartheta^i_1 - \vartheta^* - (\vartheta^i_2 - \vartheta^*) \big)
+ \tau^2 \int_{\Omega} |\nabla \vartheta^i|^2
\\ \label{nuoobvauni1}
&+ \| \chi^i \|^2_H
+ \tau \int_{\Omega} |\nabla \chi^i|^2
+ \tau ( \xi^i_1 - \xi^i_2, \chi^i_1 - \chi^i_2)
= {}- \tau \big( \pi(\chi^i_1) - \pi(\chi^i_2),  \chi^i_1 - \chi^i_2 \big).
\end{align} 
Since $\ln$, $A$ and $\beta$ are monotone, in view of \eqref{D-ini-3}
and \eqref{D-ini-4}
the second, the third and the seventh term on the left-hand side of 
\eqref{nuoobvauni1} are nonnegative. \pier{Besides}, if $\tau \leq 1/ (2C_{\pi})$, thanks to the Lipschitz continuity 
of $\pi$, the right-hand side of \eqref{nuoobvauni1} can be estimated as 
\begin{equation} \label{nuoobvauni2}
\colli{{}- \tau \big( \pi(\chi^i_1) - \pi(\chi^i_2),  \chi^i_1 - \chi^i_2 \big)} 
\leq \frac{1}{2}  \| \chi^i \|^2_H.
\end{equation} 
Then, due to \eqref{nuoobvauni2}, from \eqref{nuoobvauni1} we infer that 
\begin{equation} 
\tau^{3/2} \| \vartheta^i \|^2_H 
+ \tau^2 \int_{\Omega} |\nabla \vartheta^i|^2
+ \frac{1}{2}\| \chi^i \|^2_H
+ \tau \int_{\Omega} |\nabla \chi^i|^2
\leq 0,
\end{equation} 
whence we easily conclude that $\vartheta^i=\chi^i = 0$, i.e., $\vartheta^i_1 = \vartheta^i_2$
and $\chi^i_1 = \chi^i_2$ \colli{for $i=1,\ldots, N$.}

\setcounter{equation}{0}
\section{A priori estimates on $(AP_{\tau})$}
In this section we deduce some \pier{uniform estimates, independent of $\tau$ and}  
inferred from the equations \eqref{D-ini-0}--\eqref{D-ini-6} of 
the approximating problem $(P_{\tau})$.

\paragraph{First \pier{uniform} estimate.}
We test \eqref{D-ini-1} by $\vartheta^i$ and \eqref{D-ini-2} by $(\chi^i - \chi^{i-1}) / 
\tau$, then we sum \colli{up}. Adding $(\chi^i, \chi^i - \chi^{i-1})$ to both \colli{sides} of the resulting 
equality and exploiting the cancellation of the suitable corresponding terms, we obtain 
that
\begin{align} \nonumber
&\tau^{1/2} (\vartheta^i -  \vartheta^{i-1}, \vartheta^i)
+ \big( \ln{\vartheta^i} - \ln{\vartheta^{i-1}}, \vartheta^i \big)
+ \tau (\zeta^i, \vartheta^i - \vartheta^*)
+ \tau k_0 \int_{\Omega} | \nabla \vartheta^i |^2
\\ \nonumber
&+ \tau \bigg\| \frac{\chi^i - \chi^{i-1}}{\tau}\bigg\|^2_H
+ (\chi^i, \chi^i - \chi^{i-1})
+ (\nabla \chi^i, \nabla \chi^i - \nabla \chi^{i-1} )
+ (\xi^i,  \chi^i - \chi^{i-1} )
\\ \label{stimacomple1}
&= - \tau (\zeta^i, \vartheta^*) 
+ \tau (F^i, \vartheta^i) 
- \big(\pi(\chi^i) - \chi^i , \chi^i - \chi^{i-1} \big).
\end{align}
Due to \eqref{dis777}, we can rewrite the first, the fifth and the sixth term on the left-hand side of \eqref{stimacomple1} as
\begin{gather} \label{help1}
\tau^{1/2} (\vartheta^i -  \vartheta^{i-1}, \vartheta^i) 
= \frac{\tau^{1/2}}{2} \| \vartheta^i \|^2_H - \frac{\tau^{1/2}}{2} \| \vartheta^{i-1} \|^2_H 
+ \frac{\tau^{1/2}}{2} \| \vartheta^i - \vartheta^{i-1} \|^2_H,
\\
(\chi^i, \chi^i - \chi^{i-1})
+ (\nabla \chi^i, \nabla \chi^i - \nabla \chi^{i-1} ) = \frac{1}{2} \|\chi^i \|_V^2
- \frac{1}{2} \| \chi^{i-1} \|_V^2 + \frac{1}{2} \| \chi^i-\chi^{i-1} \|_V^2.
\end{gather}
Moreover, since the function $u \longmapsto e^u$ is convex and \pier{$e^u$ turns out to be  its subdifferential, by setting 
$u^i = \ln{\vartheta^i} $} we obtain that
\begin{equation} 
\big( \ln{\vartheta^i} - \ln{\vartheta^{i-1}}, \vartheta^i \big)
= \big( u^i - u^{i-1}, e^{u^i} \big)
\geq \int_{\Omega} e^{u^i} - \int_{\Omega} e^{u^{i-1}} 
= \| \vartheta^i \|_{L^1(\Omega)} - \| \vartheta^{i-1} \|_{L^1(\Omega)}. 
\end{equation}
Recalling that $A$ is a maximal monotone operator and $0 \in A(0)$, \pier{by \eqref{D-ini-3}
the third term on the left-hand side of \eqref{stimacomple1} is nonnegative.
We also notice that, since $\beta$ is the subdifferential of $\tilde{\beta}$, from 
\eqref{D-ini-4} it follows that}
\begin{equation} 
(\xi^i,  \chi^i - \chi^{i-1} ) \geq \int_{\Omega} \tilde{\beta}(\chi^i) - \int_{\Omega} \tilde{\beta}(\chi^{i-1}),
\end{equation}
while, due to \eqref{Poincare}, \eqref{dis2}
and the sub-linear growth of $A$ stated by \eqref{stimaA}, we deduce that
\begin{gather} \nonumber
- \tau (\zeta^i,  \vartheta^* ) \leq  C_A \tau(1 + \|\vartheta^i -  \vartheta^* \|_H) \|\vartheta^* \|_H
\leq c \tau (1 + \| \vartheta^i \|_H )
\leq c \tau (1 + \| \vartheta^i \|_V )
\\ 
\leq c \tau (1 + \| \vartheta^i \|_{L^1(\Omega)} +  \| \nabla \vartheta^i \|_H)
\leq c \tau + \tau C_1 \| \vartheta^i \|_{L^1(\Omega)} +  \tau \frac{k_0}{2}  \| \nabla \vartheta^i \|^2_H,
\end{gather}
where \pier{we have applied the Young inequality in the last term and where the constant $C_1$ depends on $C_A$, $\| \vartheta^* \|_H$ and $\pier{C_p}$.}
Due to the the boundedness of $F^i$ in $L^{\infty}(\Omega)$
and the Lipschitz continuity of $\pi$, we also infer that
\begin{gather} 
\pier{\tau (F^i, \vartheta^i) \leq \tau \| F^i \|_{L^{\infty}(\Omega)} \| \vartheta^i \|_{L^1(\Omega)}} ,
\\[0.3cm] \nonumber
- \big(\pi(\chi^i) - \chi^i, \chi^i - \chi^{i-1} \big) 
\leq c \tau (1 + \| \chi^i \|_H) \bigg\| \frac{\chi^i - \chi^{i-1}}{\tau}\bigg\|_H \\ \label{help2}
\leq \frac{\tau}{2}\bigg\| \frac{\chi^i - \chi^{i-1}}{\tau}\bigg\|^2_H +  \tau C_2 (1 + \| \chi^i \|^2_H),
\end{gather}
where $C_2$ depends on \colli{$C_{\pi}$, $|\pi(0)|$ and $|\Omega|$}.
Now, we apply the estimates \eqref{help1}--\eqref{help2} to the corresponding terms of \eqref{stimacomple1} and
sum up for $i= 1, \colli{\ldots\,}, n$, letting $n \leq N$. We obtain that
\begin{align} \nonumber
&\frac{\tau^{1/2}}{2} \| \vartheta^n \|^2_H 
+  \sum^n_{i=1} \frac{\tau^{1/2}}{2} \| \vartheta^i -  \vartheta^{i-1} \|^2_H
+ \| \vartheta^n \|_{L^1(\Omega)}
+ \frac{k_0}{2} \sum^n_{i=1} \tau \| \nabla \vartheta^i \|^2_H
\\ \nonumber
&{}+ \frac{1}{2} \sum^n_{i=1} \tau \bigg\| \frac{\chi^i - \chi^{i-1}}{\tau}\bigg\|^2_H
+ \frac{1}{2} \| \chi^n \|^2_V
+ \frac{1}{2} \sum^n_{i=1} \| \chi^i - \chi^{i-1}  \|^2_V
+ \int_{\Omega} \tilde{\beta}(\chi^n)
\\ \nonumber
&{}\leq 
\frac{\tau^{1/2}}{2} \| \vartheta_0 \|^2_H 
+ \| \vartheta_0 \|_{L^1(\Omega)}
+ \frac{1}{2} \| \chi_0 \|^2_V
+ \int_{\Omega} \tilde{\beta}(\chi_0)
+ \tau \sum^n_{i=1}  \| F^i \|_{L^{\infty}(\Omega)} \| \vartheta^i \|_{L^1(\Omega)}
\\ \label{stimacomple2}
&\quad{}+ C_1 \sum^n_{i=1} \tau \| \vartheta^i \|_{L^1(\Omega)}
+ C_2 \sum^n_{i=1} \tau\| \chi^i \|^2_H + c.
\end{align}
\pier{On account of \eqref{teta0}--\eqref{chi0} and \eqref{chi0c}}, the first four terms on the right-hand side of  \eqref{stimacomple2}
are bounded. Now, recalling the definition \eqref{Fapp-i} of $F^i$, we have that
\begin{equation} \nonumber
\tau \sum^{n}_{i=1}  \| F^i \|_{L^{\infty}(\Omega)}  \| \vartheta^i \|_{L^1(\Omega)}
=  \| \vartheta^n \|_{L^1(\Omega)} \int_{(n-1) \tau}^{n \tau} \| F(s) \|_{L^{\infty}(\Omega)} \ ds  
+ \sum^{n-1}_{i=1} \| F^i \|_{L^{\infty}(\Omega)}  \| \vartheta^i \|_{L^1(\Omega)}.
\end{equation}
Thanks to the absolute continuity of the integral, if \pier{$\tau$ is small enough (independently of $n$)} we have that
\begin{equation} \label{pier7}
\int_{(n-1) \tau}^{n \tau} \| F(s) \|_{L^{\infty}(\Omega)} \ ds 	\leq  \frac{1}{4}, 
\quad \  C_1 \tau \leq  \frac{1}{4},
\quad \  C_2 \tau \leq  \frac{1}{4}.
\end{equation}
Then, on the basis of \eqref{pier7}, from  \eqref{stimacomple2} we infer that
\begin{align} \nonumber
&\frac{\tau^{1/2}}{2} \| \vartheta^n \|^2_H 
+  \sum^n_{i=1} \frac{\tau^{1/2}}{2} \| \vartheta^i -  \vartheta^{i-1} \|^2_H
+ \frac{1}{2} \| \vartheta^n \|_{L^1(\Omega)}
+ \frac{k_0}{2} \sum^n_{i=1} \tau \| \nabla \vartheta^i \|^2_H
\\ \nonumber
&{}+ \frac{1}{2} \sum^n_{i=1} \tau \bigg\| \frac{\chi^i - \chi^{i-1}}{\tau}\bigg\|^2_H
+ \frac{1}{4} \| \chi^n \|^2_V
+ \frac{1}{2} \sum^n_{i=1} \| \chi^i - \chi^{i-1}  \|^2_V
+ \int_{\Omega} \tilde{\beta}(\chi^n)
\\
&{}\leq c  + \sum^{n-1}_{i=1} \tau \pier{\bigg( \| F^i \|_{L^{\infty}(\Omega)}  \| \vartheta^i \|_{L^1(\Omega)} + C_1 \| \vartheta^i \|_{L^1(\Omega)}
+ C_2 \| \chi^i \|^2_H \bigg).}
\end{align}
Now, we observe that
\begin{equation*} 
\sum^{n-1}_{i=1} \tau  C_1 \leq \sum^{N}_{i=1} \tau  C_1 =   C_1 T, \quad \quad
\sum^{n-1}_{i=1} \tau  C_2 \leq \sum^{N}_{i=1} \tau  C_2 =   C_2 T.
\end{equation*}
and, according to \eqref{R}, 
\begin{equation*} 
\sum^{n-1}_{i=1} \tau \| F^i \|_{L^{\infty}(\Omega)} 
\leq \sum^{N}_{i=1} \int_{(i-1) \tau}^{i \tau } \| \colli{F(s)} \|_{L^{\infty}(\Omega)} \ ds
= \int_0^T \| \colli{F (s) }\|_{L^{\infty}(\Omega)} \ ds \leq c.
\end{equation*}
Then, we \pier{can apply Lemma~\ref{grondiscr} and, recalling the notations \eqref{regtratti1},} we conclude that
\begin{align} \nonumber
&\tau^{1/2} \| \overline{\vartheta}_{\tau} \|^2_{L^{\infty}(0,T; H)} 
+ \tau^{3/2} \| \partial_t \widehat{\vartheta}_{\tau} \|^2_{L^2(0,T; H)} 
+ \| \overline{\vartheta}_{\tau} \|_{L^{\infty}(0,T; L^1(\Omega))} 
+ \colli{\| \nabla \overline{\vartheta}_{\tau} \|_{L^2(0,T; H)}^2}
\\ \label{gaf1}
&{}+ \| \partial_t \widehat{\chi}_{\tau} \|^2_{L^2(0,T; H)}
+ \| \overline{\chi}_{\tau} \|^2_{L^{\infty}(0,T; V)} 
+ \tau \| \partial_t \widehat{\chi}_{\tau} \|^2_{L^2(0,T; V)}
+ \colli{\| \tilde{\beta}(\overline{\chi}_{\tau}) \|_{L^{\infty}(0,T; L^1(\Omega))}}
\leq c. 
\end{align}
Since the third and the fourth term of the left-hand side of \eqref{gaf1} are bounded, owing to \eqref{Poincare} we also infer that
\begin{equation} \label{gaf2}
\| \overline{\vartheta}_{\tau} \|_{L^2(0,T; V)}  \leq c.
\end{equation}
\pier{Besides, in view of \eqref{D-ini-3} and due to the sub-linear growth of $A$ stated by \eqref{stimaA} and to \eqref{kappa}, we deduce that}
\begin{equation} \label{gaf3}
\| \overline{\zeta}_{\tau} \|_{L^2(0,T; H)}  \leq c.
\end{equation}

\paragraph{Second \pier{uniform} estimate.}
We formally test \eqref{D-ini-2} by $\xi^i$ and obtain
\begin{equation} \label{formale1}
(\chi^i - \chi^{i-1}, \xi^i ) + \tau \| \xi^i \|^2_H \leq \tau \big( \pi(\chi^i) + \ell \vartheta^i, \xi^i \big).
\end{equation}
We point out that \pier{the previous estimate \eqref{formale1} can be 
rigorously derived by} testing \eqref{E-D-ini-2} by 
$\beta_{\varepsilon}(X_{\varepsilon})$
and then passing to the limit as $\varepsilon \searrow 0$.
Since $\beta$ is the subdifferential of $\tilde{\beta}$, we have that
\begin{equation} \label{formale2}
(\chi^i - \chi^{i-1}, \xi^i ) \geq \int_{\Omega} \tilde{\beta}(\chi^i) - \int_{\Omega} \tilde{\beta}(\chi^{i-1}).
\end{equation}
Moreover, due to the Lipschitz continuity of $\pi$,
applying the Young inequality \eqref{dis2} to the right-hand side of \eqref{formale1}, we deduce that
\begin{equation} \label{formale3}
\tau \big( \pi(\chi^i) + \ell \vartheta^i, \xi^i \big) \leq \frac{1}{2} \tau \| \xi^i \|^2_H
+ c \, \tau \big( 1 + \| \chi^i \|^2_H + \| \vartheta^i \|^2_H \big).
\end{equation}
Now, combining \eqref{formale1}--\eqref{formale3} and summing up for $i=1, \colli{\ldots\,}, n$, with $n\leq N$, we infer~that
\begin{equation} \label{formale4}
\int_{\Omega} \tilde{\beta}(\chi^n) 
+ \frac{1}{2} \sum^n_{i=1} \tau \| \xi^i \|^2_H
\leq \int_{\Omega} \tilde{\beta}(\chi_0)
+ \sum^n_{i=1} \tau \big( 1 + \| \chi^i \|^2_H + \| \vartheta^i \|^2_H \big),
\end{equation}
whence, due to \eqref{gaf1}--\eqref{gaf2}, we obtain that
\begin{equation} \label{formale5}
\| \overline{\xi}_{\tau} \|_{L^2(0,T; H)} \leq c.
\end{equation}
Finally, by comparison in \eqref{appr2}, we conclude that $\| \Delta 
\overline{\chi}_{\tau} \|_{L^2(0,T; H)} \leq c$. \pier{Then, thanks to 
\eqref{gaf1} and elliptic regularity, we find that}
\begin{equation} \label{formale7}
\| \overline{\chi}_{\tau} \|_{L^2(0,T; W)} \leq c.
\end{equation}

\paragraph{Third \pier{uniform} estimate.}
We \pier{introduce} the function $\psi_n: \mathbb{R} \longmapsto \mathbb{R}$ 
obtained by truncating \pier{the logarithmic function} in the following way:
\begin{displaymath}
\psi_n (u) =
\left\{
\begin{array}{ll}
\ln (u)
&\qquad\text{if $u \geq 1/n$,}  \\
\noalign{\medskip}
- \ln (n)
&\qquad\text{if $u < 1/n$}.
\end{array}
\right.
\end{displaymath}
\pier{It is easy to see that $\psi_n$ is an increasing and} Lipschitz 
continuous function. Then, defining 
\begin{equation} \label{j}
j_n (u) = \int_1^u \psi_n (s) \ ds, \quad \  u \in \mathbb{R}, \quad \  \textrm{and} \  \quad
j (u) = \int_1^u \ln s \ ds, \quad \  u > 0,
\end{equation}
and testing  \eqref{D-ini-1} by $\psi_n (\vartheta^i)$, we obtain that
\begin{align}\nonumber
&\tau^{1/2} \big(\vartheta^i - \vartheta^{i-1}, \psi_n (\vartheta^i)\big)
+ \big( \ln{\vartheta^i} - \ln{\vartheta^{i-1}} , \psi_n (\vartheta^i)\big)
+ \tau k_0 \int_{\Omega \cap \{ \vartheta^i \geq 1/n \}} \frac{| \nabla \vartheta^i |^2}{\vartheta^i} 
\\ \label{sb10}
&{}= - \ell \big(\chi^i - \chi^{i-1}, \psi_n (\vartheta^i)\big)
- \tau \big(\zeta^i , \psi_n (\vartheta^i)\big) 
+ \tau \big(F^i , \psi_n (\vartheta^i)\big).
\end{align}
Recalling that \pier{$j_n$ is a convex function with derivative  $\psi_n$,} we have that
\begin{equation*}
\tau^{1/2} \big(\vartheta^i - \vartheta^{i-1}, \psi_n (\vartheta^i)\big)
\geq \tau^{1/2} \int_{\Omega} j_n( \vartheta^i) - \tau^{1/2} \int_{\Omega} j_n(\vartheta^{i-1}),
\end{equation*}
and consequently from \eqref{sb10} we infer that
\begin{align}\nonumber
&\tau k_0 \int_{\Omega \cap \{ \vartheta^i \geq 1/n \}} \frac{| \nabla \vartheta^i |^2}{\vartheta^i} 
\leq \tau^{1/2} \int_{\Omega} j_n( \vartheta^{i-1}) - \tau^{1/2} \int_{\Omega} j_n(\vartheta^i)
\\ \label{sb12}
&{}- \int_{\Omega} ( \ln{\vartheta^i} - \ln{\vartheta^{i-1}}) \psi_n (\vartheta^i)
- \int_{\Omega} \big( \ell(\chi^i - \chi^{i-1}) + \tau \zeta^i - \tau F^i \big) \psi_n (\vartheta^i).
\end{align}
Due to the properties of the subdifferential, we have that
\begin{equation} \label{pier8}
\pier{0 \leq j( \vartheta^k) \leq j(1) + (\ln \vartheta^k, \vartheta^k - 1) \quad \hbox{ for } \,
k=0,1, \ldots, N.}
\end{equation}
Since \pier{$\ln{\vartheta^k} \in H$, $\vartheta^k > 0$ a.e. in $\Omega$ and $ \vartheta^k \in H$, 
from \eqref{pier8} we infer that $j(\vartheta^k) \in L^1(\Omega)$; consequently, passing to the limit as $ n \rightarrow + \infty$, we obtain that
\begin{eqnarray*}
\psi_n (\vartheta^k) & \rightarrow & \ln{\vartheta^k} \quad \quad \textrm{ in $H$ and a.e. in $\Omega$,}  \\
j_n( \vartheta^k)    & \rightarrow & j(\vartheta^k) \quad \quad \textrm{in $L^1(\Omega)$ and a.e. in $\Omega$,}  
\end{eqnarray*}
for $
k=0,1, \ldots, N.$}
Then, taking the $\liminf$ in \eqref{sb12} as $n \rightarrow + \infty$ 
and applying the Fatou Lemma and \eqref{dis777}, we have that 
\begin{align}\nonumber
&\tau k_0 \int_{\Omega} \frac{| \nabla \vartheta^i |^2}{\vartheta^i} 
\leq \tau^{1/2} \int_{\Omega} j( \vartheta^{i-1}) - \tau^{1/2} \int_{\Omega} j(\vartheta^i)
+ \frac{1}{2} \int_{\Omega} | \ln{\vartheta^{i-1}} |^2
- \frac{1}{2} \int_{\Omega} | \ln{\vartheta^{i}} |^2
\\ \label{sb13}
&{}- \frac{1}{2} \int_{\Omega} | \ln{\vartheta^i} - \ln{\vartheta^{i-1}} |^2- \int_{\Omega} \big( \ell(\chi^i - \chi^{i-1}) + \tau \zeta^i - \tau F^i \big) \ln{\vartheta^i}.
\end{align}
Now, sum up \eqref{sb13} for $i= 1, \colli{\ldots\,}, k$, with $k \leq N$, and obtain that
\begin{align}\nonumber
&\tau^{1/2} \int_{\Omega} j( \vartheta^k) 
+ \frac{1}{2} \| \ln{\vartheta^k} \|_H^2
+ \frac{1}{2} \sum_{i=1}^k \tau^2  \bigg\|  \frac{ \ln{\vartheta^i} - \ln{\vartheta^{i-1}} }{\tau} \bigg\|_H^2
+  k_0 \sum_{i=1}^k \tau \int_{\Omega} \frac{| \nabla \vartheta^i |^2}{\vartheta^i} 
\\ \nonumber
&{}\leq \tau^{1/2} \int_{\Omega} j( \vartheta_0) 
+ \frac{1}{2} \| \ln{\vartheta_0} \|_H^2
+ \frac{1}{4} \sum_{i=1}^k \tau \| \ln{\vartheta^i } \|_H^2
+ c \sum_{i=1}^k \tau \bigg\| \frac{\chi^i - \chi^{i-1}}{\tau} \bigg\|^2_H
\\ \label{sb14}
&\quad{} + c \sum_{i=1}^k \tau \| \zeta^i \|^2_H
+ c \sum_{i=1}^k \tau \| F^i \|^2_H.
\end{align}
We observe that if $\tau \leq 1$ then
\begin{align} \frac{1}{4} \sum_{i=1}^k \tau \| \ln{\vartheta^i } \|_H^2 
\label{sb15}
\leq \frac{1}{4} \sum_{i=1}^{k-1} \tau \| \ln{\vartheta^i } \|_H^2 +
\frac{1}{4}  \| \ln{\vartheta^k } \|_H^2 \colli{.}
\end{align}
We also notice that the fourth and the fifth term on the right-hand side of \eqref{sb14}
are bounded by a positive constant $c$, due to \eqref{gaf1} and \eqref{gaf3}, respectively.
Moreover, thanks to \eqref{R} and to the definition \eqref{Fapp-i} of $F^i$,
by using the H\"older inequality the last term on the right-hand side of \eqref{sb14} can be estimated as follows:
\begin{align} \nonumber
&c \sum_{i=1}^k \tau \| F^i \|^2_H 
 \leq c \sum_{i=1}^k \tau \bigg\| \frac{1}{\tau} \int_{(i-1)\tau}^{i \tau} F(s) \ ds \bigg\|^2_H 
 \\  \label{sb16}
&{} \leq c \sum_{i=1}^k  \int_{(i-1)\tau}^{i \tau} \|   F(s) \|^2_H \ ds \leq  c \|   F \|^2_{L^2(0,T;H)}.
\end{align}
Then, combining \eqref{sb14} with \eqref{sb15}--\eqref{sb16}
\pier{(see also \eqref{teta0} and \eqref{pier8})}, we infer that
\begin{align*}\nonumber
&\tau^{1/2} \int_{\Omega} j( \vartheta^k) 
+ \frac{1}{4} \| \ln{\vartheta^k} \|_H^2
+ \frac{1}{2} \sum_{i=1}^k \tau^2  \bigg\|  \frac{ \ln{\vartheta^i} - \ln{\vartheta^{i-1}} }{\tau} \bigg\|_H^2 
\\  
&{}+ 4  k_0 \sum_{i=1}^k \tau \int_{\Omega} \left| \nabla (\vartheta^i)^{1/2} \right|^2\leq \pier{{}c
+ \frac{1}{4} \sum_{i=1}^{k-1} \tau \| \ln{\vartheta^i } \|_H^2},
\end{align*}
whence, by applying Lemma~\ref{grondiscr}, we conclude that
\begin{equation} \label{sb17}
\tau^{1/2} \| j( \overline{\vartheta}_{\tau}) \|_{L^{\infty}(0,T; L^1(\Omega))}
+ \| { \overline{\ln\vartheta}_{\tau} } \|_{L^{\infty}(0,T; H)}
+ \big\| \nabla \overline{\vartheta^{1/2}}_{\tau} \big\|_{L^2(0,T; H)} \leq c.
\end{equation}
Moreover, due to \eqref{gaf1} \colli{as well,} we also infer that
\begin{equation} \label{sb18}
\|  \overline{\vartheta^{1/2}}_{\tau} \|_{L^{\infty}(0,T; H) \cap L^2(0,T; V)} \leq c.
\end{equation}

\paragraph{Fourth \pier{uniform} estimate.}
We test \eqref{D-ini-1} by $(\vartheta^i - \vartheta^{i-1})$. \pier{Then, we take the difference between  
\eqref{D-ini-2} written for $i$ and for $i-1$, and test by $(\chi^i - \chi^{i-1})/ \tau$. Using \eqref{A1-bis} and adding, it is note difficult to} obtain that
\begin{align} \nonumber
&\tau^{1/2} \| \vartheta^i - \vartheta^{i-1} \|^2_H
+ ( \ln{\vartheta^i} - \ln{\vartheta^{i-1}}, \vartheta^i - \vartheta^{i-1} )
+ \ell (\chi^i - \chi^{i-1}, \vartheta^i - \vartheta^{i-1}) 
\\ \nonumber
&{}+ \tau \Phi (\vartheta^i - \vartheta^*)
- \tau \Phi (\vartheta^{i-1}- \vartheta^*)
+ \tau \frac{k_0}{2}\left( \| \nabla \vartheta^i\|^2_H
+\| \nabla (\vartheta^i - \vartheta^{i-1})\|^2_H
- \| \nabla \vartheta^{i-1} \|^2_H \right)
\\ \nonumber
&{}+ \frac{\tau}{2} \bigg\| \frac{\chi^i - \chi^{i-1}}{\tau}   \bigg\|^2_H
+ \frac{\tau}{2} \bigg\| \frac{\chi^i - \chi^{i-1}}{\tau}   - \frac{\chi^{i-1} - \chi^{i-2}}{\tau}\bigg\|^2_H
- \frac{\tau}{2} \bigg\| \frac{\chi^{i-1} - \chi^{i-2}}{\tau}\bigg\|^2_H
\\\nonumber
&{}+ \tau^2 \bigg\| \nabla \frac{\chi^i - \chi^{i-1}}{\tau} \bigg\|^2_H
+ (\xi^i - \xi^{i-1}, \chi^i - \chi^{i-1})
- \tau \ell \bigg( \vartheta^i - \vartheta^{i-1}, \frac{\chi^i - \chi^{i-1}}{\tau} \bigg)
\\ \label{sommona1}
&\leq \tau ( F^i, \vartheta^i - \vartheta^{i-1})
- \tau \bigg( \pi(\chi^{i}) - \pi(\chi^{i-1}), \frac{\chi^i - \chi^{i-1}}{\tau} \bigg),
\end{align}
for $i=2, \colli{\ldots\,}, N$. Now, we write \eqref{D-ini-1} and \eqref{D-ini-2}
for $i=1$ and test the corresponding equations by $(\vartheta^1 - \vartheta^0)$
and $(\chi^1 - \chi^0)/ \tau$, respectively.
Since $\vartheta^0 = \vartheta_0$ and $\chi^0=\chi_0$, we have that
\begin{align} \nonumber
&\tau^{1/2} \| \vartheta^1 - \vartheta^0 \|^2_H
+ ( \ln{\vartheta^1} - \ln{\vartheta^0}, \vartheta^1 - \vartheta^0 )
+ \ell (\chi^1 - \chi^0, \vartheta^1 - \vartheta^0) 
+ \tau \Phi (\vartheta^1 - \vartheta^*)
\\ \nonumber
&{}- \tau \Phi (\vartheta_0- \vartheta^*)
+ \tau \frac{k_0}{2} \left( \| \nabla \vartheta^1\|^2_H
+ \| \nabla (\vartheta^1 - \vartheta^0)\|^2_H
- \| \nabla \vartheta_0 \|^2_H \right)
+ \tau \bigg\| \frac{\chi^1 - \chi^0}{\tau}   \bigg\|^2_H
\\ \nonumber
&{}+ \| \nabla (\pier{\chi^1 - \chi^0} )\|^2_H
+ (\xi^1 - \xi_0, \chi^1 - \chi_0)
\leq - \tau \bigg( \pi(\chi^1) - \pi(\chi^0), \frac{\chi^1 - \chi^0}{\tau} \bigg)
\\ \label{sommona2}
&{}+ \tau \ell \bigg( \vartheta^1 - \vartheta^0, \frac{\chi^1 - \chi^0}{\tau} \bigg)
+ \tau ( F^1, \vartheta^1 - \vartheta^0)
+ (\ell \vartheta_0 + \Delta \chi_0 - \xi_0 - \pi(\chi_0), \chi^1 - \chi^0).
\end{align}
Then, we divide \eqref{sommona1} and \eqref{sommona2} by $\tau$
and sum up the corresponding equations for $i=1, \colli{\ldots\,}, n$, with $n \leq N$.
Since $\beta$ is maximal monotone \pier{and \eqref{D-ini-4} and \eqref{chi0b} hold}, then
the \colli{eleventh} term on the left-hand side of \eqref{sommona1} 
and the \colli{ninth} term on the left-hand side of \eqref{sommona2} are nonnegative.
Assuming $\chi^{-1}= \chi_0$, we infer that 
\begin{align} \nonumber
&\tau^{1/2} \sum^n_{i=1} \tau \bigg\| \frac{\vartheta^i - \vartheta^{i-1}}{\tau} \bigg\|^2_H
+ \sum^n_{i=1} \frac1\tau ( \ln{\vartheta^i} - \ln{\vartheta^{i-1}}, {\vartheta^i - \vartheta^{i-1}})
+ \Phi (\vartheta^n - \vartheta^*)
\\ \nonumber
&{}+ \frac{k_0}{2} \| \nabla \vartheta^n\|^2_H
+ \frac{k_0}{2} \tau \sum^n_{i=1} \tau \bigg\| \nabla \frac{\vartheta^i - \vartheta^{i-1}}{\tau} \bigg\|^2_H
+ \frac{1}{2} \bigg\| \frac{\chi^n - \chi^{n-1}}{\tau}   \bigg\|^2_H
\\ \nonumber
&{}+ \frac{1}{2} \sum^n_{i=1} \bigg\| \frac{\chi^i - \chi^{i-1}}{\tau}   - \frac{\chi^{i-1} - \chi^{i-2}}{\tau}\bigg\|^2_H
+ \frac12 \bigg\| \frac{\chi^1 - \chi^0}{\tau}   \bigg\|^2_H+ \sum^n_{i=1} \tau \bigg\| \nabla \frac{\chi^i - \chi^{i-1}}{\tau} \bigg\|^2_H
\\ \nonumber
&{}\leq \Phi (\vartheta_0- \vartheta^*)
+  \frac{k_0}{2} \| \nabla \vartheta_0 \|^2_H
+ \| \ell \vartheta_0 + \Delta \chi_0 - \xi_0 - \pi(\chi_0)\|^2_H
+ \frac{1}{4} \bigg\| \frac{\chi^1 - \chi^0}{\tau} \bigg\|^2_H
\\ \label{sommona3}
&\quad{}
+ ( F^n, \vartheta^n)
- (F^1, \vartheta_0)
- \sum^{n-1}_{i=1} (F^{i+1} - F^i, \vartheta^i )
+\sum^{n}_{i=1} C_{\pi} \tau \bigg\| \frac{\chi^i - \chi^{i-1}}{\tau} \bigg\|^2_H.
\end{align}
\pier{In view of  \eqref{lambda}--\eqref{chi0b} and noting that} $\vartheta_0 \in V$, $\chi_0 \in W$ 
and $\Phi$ has at most a quadratic growth
(see \eqref{teta0}--\eqref{chi0} and \eqref{A1-bis}),
the first three terms on the right-hand side of \eqref{sommona3} are bounded by a positive constant.
\pier{Besides}, using \eqref{dis2}, \eqref{Poincare} and the H\"older inequality and recalling \eqref{R}
and \eqref{gaf1},
the fifth and the sixth term on the right-hand side of \eqref{sommona3} can be estimated as follows:
\begin{gather} 
\nonumber
\hskip-2cm |( F^n, \vartheta^n)| \leq \| F^n \|_H \| \vartheta^n \|_H 
\leq \pier{C_p}  \| F^n \|_H \Big( \| \vartheta^n \|_{L^1(\Omega)} + \| \nabla \vartheta^n \|_H \Big)
\\
\label{xxz1}
\leq \pier{C_p}  \| F \|_{C^0([0,T];H)} \Big( \| \overline{\vartheta}_\tau \|_{L^\infty(0,T;L^1(\Omega))} 
+ \| \nabla \vartheta^n \|_H \Big)
\leq \frac{k_0}{4} \| \nabla \vartheta^n \|^2_H + c,
\\[0.2cm]
\label{xxz2}
|( F^1, \vartheta_0)| \leq \| F^1 \|_H \| \vartheta_0 \|_H 
\leq \| F \|_{C^0([0,T];H)}   \| \vartheta_0 \|_H \leq c.
\end{gather}
With the help of \eqref{dis2}, H\"older's inequality and \eqref{gaf2} we also infer that
\begin{align}
 \nonumber
&\bigg| \sum^{n-1}_{i=1} (F^{i+1} - F^i, \vartheta^i ) \bigg|
\leq   \sum^{n}_{i=2} \tau \bigg\| \frac{F^{i} - F^{i-1}}{\tau} \bigg\|_H \| \vartheta^{i-1}  \|_H
\\ \label{xxz3}
&{}\leq   \frac{1}{2} \sum^{n}_{i=2} \tau \bigg\| \frac{F^{i} - F^{i-1}}{\tau} \bigg\|^2_H 
+\frac{1}{2} \sum^{n-1}_{i=1} \tau \| \vartheta^{i}  \|^2_H
\leq   \frac{1}{2} \sum^{n}_{i=2} \tau \bigg\| \frac{F^{i} - F^{i-1}}{\tau} \bigg\|^2_H +c.
\end{align}
Recalling \eqref{R} and the definition of $F^i$ (see \eqref{Fapp-i}), \pier{we have that 
\begin{align} \nonumber
&\bigg\| \frac{F^i - F^{i-1}}{\tau} \bigg\|^2_H 
= \bigg\| \frac{1}{\tau^2} \int_{(i-1) \tau}^{i \tau} F(s) \, ds  - \frac{1}{\tau^2} \int_{(i-2) \tau}^{(i-1) \tau} 
F(s) \, ds \bigg\|^2_H 
\\ \nonumber
&= \bigg\| \frac{1}{\tau^2} \int_{(i-1) \tau}^{i \tau} \Big( F(s) - F(s-\tau) \Big)  d s  \bigg\|^2_H \\ \nonumber
& \leq  
\frac{1}{\tau^4} \left| \int_{(i-1) \tau}^{i \tau} \| F(s) - F(s-\tau) \|_H \, d s  \right|^2 
\leq \frac{1}{\tau^3}   \int_{(i-1) \tau}^{i \tau} \bigg\| \int_{s-\tau}^{s} \partial_t F(t)\, dt \bigg\|^2_H \,d s   
\\ \nonumber
&{}\leq \frac{1}{\tau^2} \int_{(i-1) \tau}^{i \tau}  \Bigg( \int_{s-\tau}^{s} \| \partial_t 
F(t) \|^2_H \, dt \Bigg) \, d s
\leq \frac1\tau \| \partial_t F \|^2_{L^2((i-2)\tau, i\tau; H)}  ,
\end{align}
so that
\begin{equation} \label{xxz4}
 \frac{1}{2} \sum^{n}_{i=2} \tau \bigg\| \frac{F^{i} - F^{i-1}}{\tau} \bigg\|^2_H
\leq  
 \| \partial_t F \|^2_{L^2(0,T; H)} .
\end{equation}
Next, we take advantage of Lemma~\ref{superlemma1} in order to deal with the second term on the left-hand side of \eqref{sommona3}. Indeed (cf.~\eqref{pier2}), we realize that
\begin{equation*} 
\left| (\vartheta^i)^{1/2} - (\vartheta^{i-1})^{1/2}\right| ^2\leq (\ln \vartheta^i - \ln \vartheta^{i-1}, \vartheta^i -\vartheta^{i-1}),
\end{equation*}
whence
\begin{equation}
\label{xxz5}
\sum^n_{i=1} \frac1\tau ( \ln{\vartheta^i} - \ln{\vartheta^{i-1}}, {\vartheta^i - \vartheta^{i-1}})
\geq \sum^n_{i=1} \tau \left\|   \frac{ (\vartheta^i)^{1/2} - (\vartheta^{i-1})^{1/2}}\tau    \right\|_H^2.
\end{equation}
Collecting now  \eqref{xxz1}--\eqref{xxz5}, from \eqref{sommona3} and \eqref{gaf1}} we infer that
\begin{align}\nonumber
&\tau^{1/4} \| \partial_t \widehat{\vartheta}_{\tau} \|_{L^2(0, T; H)}
+ \| \partial_t \widehat{\vartheta^{1/2}}_{\tau} \|_{L^2(0, T; H)}
+ \| \Phi (\overline{\vartheta}_{\tau} - \vartheta^*) \|_{L^{\infty}(0,T)}
+ \|\overline{\vartheta}_{\tau} \|_{L^{\infty}(0,T;V)}
\\ \label{aaa1}
&{}+ \pier{\tau^{1/2}}\| \partial_t \widehat{\vartheta}_{\tau} \|_{L^2(0, T; V)}
+ \| \partial_t \widehat{\chi}_{\tau} \|_{L^{\infty}(0, T; H)}
+ \| \partial_t \widehat{\chi}_{\tau} \|_{L^2(0, T; V)}
\leq c.
\end{align}
\pier{Therefore, thanks to \eqref{appr3} and using \eqref{stimaA} and \eqref{kappa}, we have that
\begin{equation} \label{pier9}
\| \overline{\zeta}_{\tau} \|_{L^\infty (0,T; H)}  \leq c.
\end{equation}
Moreover,} by comparison in \eqref{appr1} \pier{and in view of} \eqref{gaf1}--\eqref{gaf3},
\eqref{formale5}--\eqref{formale7}, \eqref{sb17}--\eqref{sb18} and \eqref{aaa1}, we obtain that
\begin{align} \nonumber
&\| \partial_t \widehat{\ln \vartheta}_{\tau} \|_{L^2(0, T; V')}
\leq c \tau^{1/2} \| \partial_t \widehat{\vartheta}_{\tau} \|_{L^2(0, T; H)}
+ c \| \partial_t \widehat{\chi}_{\tau} \|_{L^2(0, T; H)}
\\ \label{aaa2}
&\quad {}+ c \|\overline{\zeta}_{\tau} \|_{L^2(0,T;H)}
+k_0 \|\Delta \overline{\vartheta}_{\tau} \|_{L^2(0,T;V')}
+ c \|\overline{F}_{\tau} \|_{L^2(0,T;H)} \leq c.
\end{align}
Furthermore, recalling \eqref{appr2}, a comparison of the terms yields the bound  
$$ \| \Delta \overline{\chi}_{\tau} + \overline{\xi}_{\tau} \|_{L^{\infty}(0, T; H)} \leq c . $$
Hence, by arguing as in the Second uniform estimate, we can improve \eqref{formale5} and \eqref{formale7} 
to find out that
\begin{equation} \label{pier10}
\| \overline{\xi}_{\tau} \|_{L^\infty (0,T; H)} +  \| \overline{\chi}_{\tau} \|_{L^\infty  (0,T; W)} \leq c.
\end{equation}

\paragraph{Summary of the \pier{uniform} estimates.}
Let us \pier{collect the previous estimates. From 
\eqref{gaf1}--\eqref{gaf3},
\eqref{formale5}--\eqref{formale7}, \eqref{sb17}--\eqref{sb18} and \eqref{aaa1}--\eqref{pier10} 
we conclude that there exists a constant $c>0$, independent of $\tau$, such that
\begin{align} \nonumber
&\| \overline{\vartheta}_{\tau} \|_{L^{\infty}(0,T; V)} 
+ \| \widehat{\vartheta}_{\tau} \|_{L^{\infty}(0,T; V)}
+ \tau^{1/4} \| \partial_t \widehat{\vartheta}_{\tau} \|_{ L^2(0, T; H)}
\\ \nonumber
&{}+ \|  \overline{\ln\vartheta}_{\tau}  \|_{L^{\infty}(0,T; H)}
+ \| \widehat{\ln{\vartheta}}_{\tau}  \|_{H^1(0,T; V') \cap L^{\infty}(0,T; H)}
\\ \nonumber
&{}+ \|  \overline{\vartheta^{1/2}}_{\tau} \|_{L^{\infty}(0,T; H) \cap L^2(0,T; V)} 
+ \| \widehat{\vartheta^{1/2}}_{\tau} \|_{H^1(0, T; H)\cap L^2(0,T; V)}
+ \| \overline{\zeta}_{\tau} \|_{L^\infty (0,T; H)} 
\\ \label{gransommafinale}
&{}+ \| \overline{\chi}_{\tau} \|_{L^{\infty}(0,T; W)} 
+ \| \widehat{\chi}_{\tau} \|_{W^{1, \infty}(0, T; H) \cap H^1(0, T; V) \cap L^{\infty}(0,T; W)}  
+ \| \overline{\xi}_{\tau} \|_{L^{\infty}(0,T; H)} 
\leq c.
\end{align}}

\setcounter{equation}{0}
\section{Passage to the limit as $\tau \searrow 0$}
Thanks to \eqref{gransommafinale} and to the well-known weak or weak* compactness results, 
we deduce that, at least for a subsequence of $\tau \searrow 0$, there exist ten limit functions 
$\vartheta$, $\widehat{\vartheta}$, 
$\lambda$, $\widehat{\lambda}$,
$w$, $\widehat{w}$, $\zeta$\colli{,}
$\chi$, $\widehat{\chi}$, and 
$\xi$  such that
\begin{eqnarray} \label{ubi1}
\overline{\vartheta}_{\tau}  \rightharpoonup^* \vartheta & \textrm{in}& L^{\infty}(0,T; V), \\ 
\label{ubi2}
\widehat{\vartheta}_{\tau}   \rightharpoonup^* \widehat{\vartheta} & \textrm{in}& L^{\infty}(0,T; V), \\ 
\label{pier11}
\tau^{1/4} \widehat{\vartheta}_{\tau}   \rightharpoonup^*   0 &\textrm{in}& H^1(0,T; H)\cap L^{\infty}(0,T; V), \\ \label{ubi3}
\overline{\ln \vartheta}_{\tau}  \rightharpoonup^* \lambda & \textrm{in}& L^{\infty}(0,T; H), \\ 
\label{pier12}
 \widehat{\ln \vartheta}_{\tau}   \rightharpoonup^* \widehat{\lambda} & \textrm{in}& H^1(0,T; V') \cap L^{\infty}(0,T; H), \\ 
\overline{\vartheta^{1/2}}_{\tau}  \rightharpoonup^* w & \textrm{in}& L^{\infty}(0,T; H) \cap L^2(0,T; V), 
\\ \label{pier13}
\widehat{\vartheta^{1/2}}_{\tau}   \rightharpoonup \widehat{w} & \textrm{in}& H^1(0,T; H) \cap L^2(0,T; V), 
\\ \label{pier14}
\overline{\zeta}_{\tau}  \rightharpoonup^* \zeta & \textrm{in}& L^{\infty}(0,T; H)\label{p14bis},\\
\overline{\chi}_{\tau}  \rightharpoonup^* \chi & \textrm{in}& L^{\infty}(0,T; W), \\ 
\label{pier15}
\widehat{\chi}_{\tau}   \rightharpoonup^* \widehat{\chi} & \textrm{in}& W^{1, \infty}(0,T; H) \cap H^1(0,T; V) \cap L^{\infty}(0,T; W), \\\label{ubi4}
\overline{\xi}_{\tau}  \rightharpoonup^* \xi & \textrm{in}& L^{\infty}(0,T; H).
\end{eqnarray} 
First, we observe that $\vartheta = \widehat{\vartheta}$:
indeed, thanks to \eqref{regtratti6} and \eqref{pier11}, we have that
\begin{equation*} 
\| \overline{\vartheta}_{\tau}  - \widehat{\vartheta}_{\tau} \|_{L^2(0,T; H)}
\leq \frac{\tau}{\sqrt{3}} \| \partial_t \widehat{\vartheta}_{\tau} \|_{L^2(0,T; H)}
\leq c \tau^{3/4}
\end{equation*}
and consequently $\overline{\vartheta}_{\tau} -  \widehat{\vartheta}_{\tau} \to 0$ strongly in 
$L^2(0,T; H)$.
Moreover, it turns out that $\lambda= \widehat{\lambda}$: in fact, on account of  \eqref{regtratti7} and \eqref{pier12} we have that
\begin{equation*} 
\|  \overline{ \ln \vartheta}_{\tau} -  \widehat{\ln \vartheta}_{\tau} \|_{L^{\infty}(0,T; V')}
\leq \tau^{1/2} \| \partial_t \widehat{\ln \vartheta}_{\tau} \|_{L^2(0,T; V')} \leq c \tau^{1/2} \colli{,}
\end{equation*}
whence
\begin{equation} \label{pier16}
\lim_{\tau \searrow 0} \| \overline{\ln \vartheta}_{\tau} -  \widehat{\ln \vartheta}_{\tau} \|_{L^{\infty}(0,T; V')} = 0.
\end{equation}
Similarly, thanks to  \eqref{regtratti6} and \eqref{pier15}, we see that
\begin{equation*} 
\| \overline{\vartheta^{1/2}}_{\tau} -  \widehat{\vartheta^{1/2}}_{\tau}  \|_{L^2(0,T; H)}
\leq \frac{\tau}{\sqrt{3}} \| \partial_t  \widehat{\vartheta^{1/2}}_{\tau} \|_{L^2(0,T; H)}
\leq c \tau,
\end{equation*}
which entails
\begin{equation}\label{pier17}
\lim_{\tau \searrow 0}  \| \overline{\vartheta^{1/2}}_{\tau}   -  \widehat{\vartheta^{1/2}}_{\tau}  \|_{L^2(0,T; H)} =0
\end{equation}
and $w= \widehat{w}$.  
Finally, we check that $\chi= \widehat{\chi}$. In the light of \eqref{regtratti7}, we have that
\begin{equation*} 
\| \overline{\chi}_{\tau} -  \widehat{\chi}_{\tau}  \|_{L^{\infty}(0,T; V)}
\leq \tau \| \partial_t  \widehat{\chi}_{\tau} \|_{L^2(0,T; V)}
\leq c \tau
\end{equation*}
and consequently 
\begin{equation}\label{pier18}
\lim_{\tau \searrow 0}  \| \overline{\chi}_{\tau} -  \widehat{\chi}_{\tau}  \|_{L^{\infty}(0,T; V)} =0.
\end{equation}
Next, in view of the convergences in \eqref{pier12},  \eqref{pier13}, \eqref{pier15} and 
owing to the strong compactness lemma stated in \cite[Lemma 8, p. 84]{Simon}, we have that
\begin{eqnarray}
\label{str-c-1}
\widehat{\ln \vartheta}_{\tau}   \rightarrow  \lambda  & \textrm{in}& C^0([0,T]; V'), \\ 
\label{str-c-2}
\widehat{\vartheta^{1/2}}_{\tau}   \rightarrow w& \textrm{in}& L^2(0,T; H) , \\  
\label{str-c-3}
\widehat{\chi}_{\tau}   \rightarrow \chi  & \textrm{in}& C^0([0,T]; V).
\end{eqnarray}
Then, \colli{by} \eqref{pier16}--\eqref{pier18} we can also conclude that 
\begin{eqnarray}
\label{pier21}
\overline{\ln \vartheta}_{\tau}   \rightarrow  \lambda  & \textrm{in}& L^\infty (0,T; V'), \\ 
\label{pier22}
\overline{\vartheta^{1/2}}_{\tau}   \rightarrow w& \textrm{in}& L^2(0,T; H) , \\  
\label{pier23}
\overline{\chi}_{\tau}   \rightarrow \chi  & \textrm{in}& L^\infty (0,T; V).
\end{eqnarray}
Thanks to \eqref{pier23} and to the \colli{Lipschitz} continuity of $\pi$, we have that 
$$ \pi(\overline{\chi}_{\tau}  )  \rightarrow \pi(\chi) \quad  \textrm{ in } \,  L^\infty (0,T; H). $$
Now, we check that $\lambda = \ln \theta $: in fact,  due to the weak convergence of $\overline{\vartheta}_{\tau}$ ensured by \eqref{ubi1} and to
the strong convergence of $\ln(\overline{\vartheta}_{\tau})$ in \eqref{pier21} (see~\eqref{ubi3} as well), we have that
\begin{equation*} 
\limsup_{\tau \searrow 0} \int_0^T\!\!\int_{\Omega} \left(\ln{\overline{\vartheta}_{\tau}} \right) \overline{\vartheta}_{\tau}
= \lim_{\tau \searrow 0} \int_0^T \langle \ln{\overline{\vartheta}_{\tau}} , \overline{\vartheta}_{\tau} \rangle 
=  \int_0^T \langle \lambda , \vartheta \rangle = \int_0^T \!\! \int_{\Omega} \lambda  \vartheta,
\end{equation*}
so that a standard tool for maximal monotone operators (cf., e.g., \cite[Lemma~1.3, p.~42]{Barbu}) ensure that
$ \lambda = \ln \vartheta $. In the light of \eqref{appr4} and of the convergences \eqref{ubi4} and \eqref{pier23}, it is even simpler to check that $\xi$ and $\chi$ satisfy \eqref{iniziale4}.

At this point, recalling also \eqref{ubi3}, \eqref{pier12}, \eqref{pier15} and passing to the limit in \eqref{appr1} and \eqref{appr2}, we arrive at \eqref{iniziale1} and \eqref{iniziale2}. In addition, note that \eqref{D-ini-6} implies that 
$\widehat{\ln \vartheta}_{\tau} (0) = \ln \vartheta_0$ and $\widehat{\chi}_{\tau} (0)= \chi_0$; thus,
thanks to \eqref{pier21} and \eqref{pier23}, passing to the limit as $\tau \searrow 0$ leads to the initial conditions \eqref{iniziale6}.

It remains to show \colli{\eqref{iniziale3}}. To this aim, we point out that \eqref{pier22} implies that, possibly taking another subsequence,  $\overline{\vartheta^{1/2}}_{\tau}   \rightarrow w $ almost everywhere in $Q$. 
Then, using \eqref{ubi1} and the Egorov theorem, it is not difficult to verify that 
$$\overline{\vartheta}_{\tau}   
= \left(\overline{\vartheta^{1/2}}_{\tau}\right)^2  \rightarrow w^2 \quad  \textrm{ a.e. in $Q$ and in }\,  L^2(0,T; H) ,
$$
as well as $\vartheta= w^2$. Details of this argument can be found, for instance, in  
\cite[Exercise~4.16, part~3, p.~123]{BrezAF}. Then, as $A$ induces a natural maximal monotone operator 
on $L^2(0,T;H)$, recalling \eqref{appr3} and observing that~\colli{(cf.~\eqref{p14bis})}
\begin{equation*} 
\limsup_{\tau \searrow 0} \int_0^T  \left (\overline{\zeta}_{\tau} , \overline{\vartheta}_{\tau} -\vartheta^*  \right)_{\!H}
= \lim_{\tau \searrow 0} \int_0^T  \left (\overline{\zeta}_{\tau} , \overline{\vartheta}_{\tau} -\vartheta^*  \right)_{\!H}
= \int_0^T  \left(\zeta, \vartheta -  \vartheta^* \right)_H , 
\end{equation*}
we easily recover \colli{\eqref{iniziale3}}. Therefore, Theorem~\ref{Teo-esistenza} is completely proved.

\end{document}